\documentclass[twoise]{article}
\usepackage{amsfonts}
\usepackage{amsmath}
\usepackage{amssymb}
\usepackage{amsthm}   %-------proof
\usepackage{empheq}   %-------to use \empheq
\usepackage{titlesec} %-------to use \newpagestyle
\allowdisplaybreaks% 允许一个公式环境中的公式自动分页
\newtheorem{defi}{Definition}[section]
\newtheorem{thm}{Theorem}[section]
\newtheorem{lem}{Lemma}[section]
\newtheorem{rem}{Remark}[section]
\usepackage{bm}
\usepackage{arydshln}

\begin{document}

\renewcommand{\thefootnote}{\fnsymbol{footnote}}

\title{Global existence of radial solutions for general semilinear hyperbolic systems in 3D}
\author{Silu Yin\footnotemark[2]\and Yi Zhou\footnotemark[2]}

%\footnotetext[1]{Corresponding author.}
\footnotetext[2] {School of Mathematical Sciences, Fudan University, 200433, P.R. China.}
\footnotetext{Email address: yins11@fudan.edu.cn(Silu Yin), yizhou@fudan.edu.cn(Yi Zhou).}
\date{}
\maketitle
\begin{abstract}
 We study the well-posedness of radial solutions for general nonlinear hyperbolic systems in three dimensions. We give a proof of the global existence of radial solutions for general semilinear hyperbolic systems in 3D under null condition, with small scaling invariant $\dot{W}^{2,1}(\mathbb{R}^3)$ data. We obtain a bilinear estimate that is effective to the hyperbolic systems which do not have any time decay. It allows us to achieve the boundedness of the weighted BV norm of the radial solution.
\end{abstract}

\section{Introduction and Main Result}\label{section1}

We consider the following general first-order semilinear hyperbolic system
\begin{align}\label{1.1}
\frac{\partial u}{\partial t}+\sum_{i=1}^3A_i\frac{\partial u}{\partial x^i}=Q(u,u),
\end{align}
where $u=(\rho_1,\cdots,\rho_l, v_1,\cdots,v_m)^T$ is the unknown function of $(t,x^1,x^2,x^3)$, $\rho_1,\cdots,\rho_l$ are scalar unknown functions valued in $\mathbb{R}$, and $v_1,\cdots,v_m$ are vector unknown functions valued in $\mathbb{R}^3$. $A_1$, $A_2$, $A_3$ are constant matrixes and $Q(u,u)=(Q^\rho(u,u),Q^v(u,u))^T$, $Q^\rho=(Q_1^\rho,\cdots,Q_l^\rho)$, $Q^v=({\bf Q}_1^v,\cdots,{\bf Q}_m^v)^T$ with each $Q_k^\rho,{\bf Q}_p^v=(Q_p^{v1},Q_p^{v2},Q_p^{v3})$ is quadratic in $u$.

System (\ref{1.1}) is equipped with the following initial data,
\begin{align}\label{1.2}
t=0:\ \ \ \ u=u_0(x),\ \ \ \ x\in \mathbb{R}^3.
\end{align}

In order to state our result precisely, we need to make some assumptions. The first assumption is the following.

[H1]: Cauchy problem (\ref{1.1}) (\ref{1.2}) is rotational invariant.

We know that many physical systems in three dimensions like compressible Euler equations, relativistic compressible Euler equations and compressible MHD equations etc, are invariant under rotation of coordinates. The goal in this paper is to consider the global existence of radial solutions to Cauchy problem (\ref{1.1}) (\ref{1.2}), which is invariant under rotation of coordinates.

Consider the corresponding hyperbolic system in one dimension
\begin{align}\label{1.3}
\frac{\partial u}{\partial t}+A_1\frac{\partial u}{\partial x^1}=Q(u,u),
\end{align}
where $u,A_1,Q(u,u)$ are defined as before in (\ref{1.1}).

We recall the following concept.
\begin{defi}
  The system (\ref{1.3}) satisfies {\upshape null condition} if any small plane wave solution $u=u(t+\xi x^1)$ ($\xi$ is a constant in $\mathbb{R}$) to the linearized system
  \begin{align*}
\frac{\partial u}{\partial t}+A_1\frac{\partial u}{\partial x^1}=0
\end{align*}
is always a solution to the original system (\ref{1.3}).
\end{defi}

Then we assume

[H2]:  System (\ref{1.3}) satisfies null condition.

In this paper, we only need to assume a weaker null condition for the corresponding one-dimensional system instead of the original system (\ref{1.1}). Roughly speaking, the assumption [H2] makes sure that waves of the same family do not interact in the nonlinear terms.

The last but not least important assumption is presented as

[H3]: The initial data satisfies
\begin{align}\label{1.4}
\|u_0\|_{\dot{W}^{2,1}(\mathbb{R}^3)}\leq\varepsilon,
\end{align}
where $\varepsilon$ is a small constant.

Let us take a look at the invariance group of the system (\ref{1.1}). Suppose that $u(t,x)$ satisfies system (\ref{1.1}), obviously $u^\lambda (t,x)$ also satisfies system (\ref{1.1}), which is defined by
  \begin{align*}
    u^\lambda(t,x)=\lambda u(\lambda t,\lambda x),%\qquad u^\lambda_0=\lambda u_0(\lambda x),
  \end{align*}
for all $\lambda>0$. That is, we have
\begin{align}\label{1.4-5}
\frac{\partial u^\lambda}{\partial t}+\sum_{i=1}^3A_i\frac{\partial u^\lambda}{\partial x^i}=Q(u^\lambda,u^\lambda).
\end{align}
Differential (\ref{1.4-5}) with respect to $\lambda$ and then take $\lambda=1$, we have
\begin{align*}
(\frac{\partial}{\partial t}+\sum_{i=1}^3A_i\frac{\partial}{\partial x^i})(L_0+2)u=(L_0+1)Q(u,u).
\end{align*}
Here $L_0$ denotes the scaling operator that introduced by Klainerman \cite{4},
\begin{align*}
  L_0=t\partial_t+x^j\partial_j,
\end{align*}
and we use Einstein's summation convention for repeated indices. Clearly, the norm given in (\ref{1.4}) is dimensionless that is scaling invariant
\begin{equation*}
  \|u_0\|_{\dot{W}^{2,1}(\mathbb{R}^3)}=\|u^\lambda_0\|_{\dot{W}^{2,1}(\mathbb{R}^3)}.
\end{equation*}

As a primary task, we study on the coefficient matrixes under assumption [H1]. We have the following lemma, the proof will be given in section 2.
\begin{lem}\label{lem1}
If system (\ref{1.1}) satisfies {\upshape [H1]}, then we system (\ref{1.1}) should satisfy
$$A_1=\left(\begin{array}{ccc;{2pt/2pt}ccccccc}
        &  &       & b_{11}^{11} &  0 &0    & \cdots & b_{1m}^{11} & 0      & 0 \\
        \multicolumn{3}{c;{2pt/2pt}}{\raisebox{2ex}[0pt]{\huge0}}& \vdots     &&    &            & \vdots      &&    \\
       &&& b_{l1}^{11} & 0      & 0      & \cdots & b_{lm}^{11} & 0  & 0 \\ \hdashline[2pt/2pt]
       c_{11}^{11} & \cdots & c_{1l}^{11} & 0      & 0      & 0      & \cdots & 0      & 0      & 0 \\
       0      & \cdots & 0      & 0      & 0      & d_{11}^{111} & \cdots & 0      & 0      & d_{1m}^{111} \\
       0      & \cdots & 0      & 0      & -d_{11}^{111}& 0      & \cdots & 0      & -d_{1m}^{111}& 0 \\
      \vdots &        &        & \vdots    &&     &       & \vdots      &&      \\
       c_{m1}^{11} & \cdots & c_{ml}^{11} & 0      & 0      & 0      & \cdots & 0      & 0      & 0 \\
       0      & \cdots & 0      & 0      & 0      & d_{m1}^{111} & \cdots & 0      & 0      & d_{mm}^{111} \\
       0      & \cdots & 0      & 0      &  -d_{m1}^{111}& 0      & \cdots & 0      & -d_{mm}^{111}& 0\\
     \end{array}\right);$$
$$A_2=\left(\begin{array}{ccc;{2pt/2pt}ccccccc}
             &  &       & 0      & b_{11}^{11} & 0      & \cdots & 0      & b_{1m}^{11} & 0 \\
       \multicolumn{3}{c;{2pt/2pt}}{\raisebox{2ex}[0pt]{\huge0}}    & \vdots &        &        &        & \vdots &        &   \\
           &  &       & 0      & b_{l1}^{11} & 0      & \cdots & 0      & b_{lm}^{11} & 0 \\ \hdashline[2pt/2pt]
       0      & \cdots & 0      & 0      & 0      &-d_{11}^{111} & \cdots & 0      & 0      &-d_{1m}^{111} \\
       c_{11}^{11} & \cdots & c_{1l}^{11} & 0      & 0      & 0      & \cdots & 0      & 0      & 0 \\
       0      & \cdots & 0      & d_{11}^{111} & 0      & 0      & \cdots &  d_{1m}^{111}& 0      & 0 \\
       \vdots &        &        & \vdots &        &        &        & \vdots &        &   \\
       0      & \cdots & 0      & 0      & 0      &-d_{m1}^{111} & \cdots & 0      & 0      &-d_{mm}^{111} \\
       c_{m1}^{11} & \cdots & c_{ml}^{11} & 0      & 0      & 0      & \cdots & 0      & 0      & 0 \\
       0      & \cdots & 0      & d_{m1}^{111} & 0      & 0      & \cdots & d_{mm}^{111} &  0     & 0\\
     \end{array}\right);$$
$$A_3=\left(\begin{array}{ccc;{2pt/2pt}ccccccc}
             &  &       & 0      & 0      & b_{11}^{11} & \cdots & 0      & 0      & b_{1m}^{11}\\
       \multicolumn{3}{c;{2pt/2pt}}{\raisebox{2ex}[0pt]{\huge0}}        & \vdots &        &        &        & \vdots &        &   \\
           &  &       & 0      & 0      & b_{l1}^{11} & \cdots & 0      & 0      & b_{lm}^{11} \\ \hdashline[2pt/2pt]
       0      & \cdots & 0      & 0      & d_{11}^{111} & 0      & \cdots & 0      &  d_{1m}^{111}&   0    \\
       0      & \cdots & 0      &-d_{11}^{111}& 0      & 0      & \cdots & -d_{1m}^{111}& 0      & 0 \\
       c_{11}^{11} & \cdots & c_{1l}^{11} & 0      & 0      & 0      & \cdots & 0      & 0      & 0 \\
       \vdots &        &        & \vdots &        &        &        & \vdots &        &   \\
       0      & \cdots & 0      & 0      & d_{m1}^{111} & 0      & \cdots & 0      & d_{mm}^{111} & 0  \\
       0      & \cdots & 0      &-d_{m1}^{111} & 0      & 0      & \cdots &-d_{mm}^{111} &  0     & 0\\
       c_{m1}^{11} & \cdots & c_{ml}^{11} & 0      & 0      & 0      & \cdots & 0      & 0      & 0 \\
 \end{array}\right),$$
 and
 \begin{align*}
Q_k^\rho(u,u)&=\sum_{i,j=1}^l\Gamma_{ij}^k\rho_i\rho_j+\sum_{s,q=1}^m \Omega_{sq}^{k11}v_s^T v_q,\\
Q_p^{v}(u,u)&=\sum_{1\leq j\leq l \atop 1\leq q\leq m} \bar{\Upsilon}_{jq}^{p11}\rho_j v_q+\sum_{s,q=1}^m \bar{\Omega}_{sq}^{p111}
v_s\times v_q,
\end{align*}
$k=1,\cdots,l$, $p=1,\cdots,m$.
\end{lem}

The main result is given as follows.
\begin{thm} \label{thm1}
Under assumptions {\upshape[H1-H3]}, the Cauchy problem (\ref{1.1}) (\ref{1.2}) has a unique global radial solution
\begin{align}\label{1.5}
\rho(t,x)=(\rho_1(t,r),\cdots,\rho_l(t,r))^T,\ v(t,x)=(\frac{x}{r}w_1(t,r),\cdots,\frac{x}{r}w_m(t,r))^T,
\end{align}
where $x=(x^1,x^2,x^3)$, $r=\sqrt{(x^1)^2+(x^2)^2+(x^3)^2}$.
\end{thm}

We discover that some linear combinations of the unknown functions $\tilde{\rho}_i$, $\tilde{w}_i$ should satisfy wave equations, with nonlinear terms like $r\partial_r[(\tilde{\rho}_a+\tilde{w}_a)(\tilde{\rho}_b+\tilde{w}_b)]$, $r\partial_t[(\tilde{\rho}_a+\tilde{w}_a)(\tilde{\rho}_b+\tilde{w}_b)]$, and $(\tilde{\rho}_a+\tilde{w}_a)(\tilde{\rho}_b+\tilde{w}_b)$ where $a,b\in\{1,\cdots,l+m\}$. Null condition ensures that $a\neq b$ in the nonlinear terms. If we take notation $F$ to represent these nonlinear terms, then by the following linear estimate that we will present in section 4
\begin{align}
&\sum\limits_{\alpha=0,1}\Big[\|L_0^\alpha\partial_t(ru(t,\cdot))\|_{L^1(\mathbb{R}_+)}+\|L_0^\alpha\partial_r(ru(t,\cdot))\|_{L^1(\mathbb{R}_+)}\Big]\notag\\
\leq& C\|u_0\|_{\dot{W}^{2,1}(\mathbb{R}^3)}+\sum\limits_{\alpha=0,1}\|L_0^\alpha F\|_{L^1([0,t]\times \mathbb{R}_+)},
\end{align}
one should give a suitable bilinear estimate on $F$ to obtain the global existence. Once the bilinear estimate is proved, we complete the proof of Theorem \ref{thm1} by the global iterative method. It is a standard process that can be followed by \cite{9}, so we omit it here. Here and everywhere in this paper, notation $C$ reprensents for some constant independent of the initial data.
\begin{thm}\label{thm2}
Let $i=a,b$, consider the following wave equations
\begin{equation*}
\left\{\begin{aligned}
&(\partial_t^2-\lambda_i^2\partial_r^2)(r\tilde{\rho}_i)=F_i(t,r)\\
&t=0: r\tilde{\rho}_i=r\tilde{\rho}_{i0}(r),\ (r\tilde{\rho}_i)_t=-\lambda_i(r\tilde{\tilde{w}}_{i0})_{rr}(r),
\end{aligned}\right.
\end{equation*}
\begin{equation*}
\left\{\begin{aligned}
&(\partial_t^2-\lambda_i^2\partial_r^2)(r\tilde{\tilde{w}}_i)=G_i(t,r)\\
&t=0: r\tilde{\tilde{w}}_i=r\tilde{\tilde{w}}_{i0}(r),\ (r\tilde{\tilde{w}}_i)_t=-\lambda_ir\tilde{\rho}_{i0}(r),
\end{aligned}\right.\end{equation*}
where $\tilde{\tilde{w}}_i$ is defined as
\begin{align*}
  \tilde{\tilde{w}}_{i}(r)=\int_0^r\tilde{w}_i(s)ds.
\end{align*}
If $\lambda_a\neq \lambda_b$ and $\lambda_a\lambda_b\neq 0$, then there exists
\begin{align}
&\sum\limits_{\alpha=0,1}\|L_0^\alpha r\partial_r[(\tilde{\rho}_a+\tilde{w}_a)(\tilde{\rho}_b+\tilde{w}_b)]\|_{L^1([0,t]\times \mathbb{R}_+)}\notag\\
\leq&C\Big(\|\tilde{\rho}_{a0}\|_{\dot{W}^{2,1}(\mathbb{R}^3)}+\|\tilde{w}_{a0}\|_{\dot{W}^{2,1}(\mathbb{R}^3)}+\sum\limits_{\alpha=0,1}\big(\|L_0^\alpha F_a\|_{L^1([0,t]\times\mathbb{R}_+)}+\|L_0^\alpha G_a\|_{L^1([0,t]\times\mathbb{R}_+)}\big)\Big)\notag\\
&\quad\cdot\Big(\|\tilde{\rho}_{b0}\|_{\dot{W}^{2,1}(\mathbb{R}^3)}+\|\tilde{w}_{b0}\|_{\dot{W}^{2,1}(\mathbb{R}^3)}+\sum\limits_{\alpha=0,1}\big(\|L_0^\alpha F_b\|_{L^1([0,t]\times\mathbb{R}_+)}+\|L_0^\alpha G_b\|_{L^1([0,t]\times\mathbb{R}_+)}\big)\Big).
\end{align}
\end{thm}

If one eigenvalue is equal to zero, without loss of generality, suppose $\lambda_a=0$, then $\tilde{\rho}_a$, $\tilde{w}_a$ satisfy ordinary differential equations
\begin{equation*}
\left\{\begin{aligned}
\tilde{\rho}_{at}&=f_{1a}\\
\tilde{w}_{at}&=f_{2a},
\end{aligned}\right.\end{equation*}
and the bilinear estimate has a little difference.
\begin{thm}\label{thm3}
  Consider
\begin{equation*}
\left\{\begin{aligned}
&\tilde{\rho}_{at}=F_a(t,r)\\
&t=0: \tilde{\rho}_a=\tilde{\rho}_{a0}(r),
\end{aligned}\right.\qquad
\left\{\begin{aligned}
&\tilde{w}_{at}=G_a(t,r)\\
&t=0:\ \tilde{w}_a=\tilde{w}_{a0}(r),
\end{aligned}\right.
\end{equation*}
\begin{equation*}
\left\{\begin{aligned}
&(\partial_t^2-\lambda_b^2\partial_r^2)(r\tilde{\rho}_b)=F_b(t,r)\\
&t=0: r\tilde{\rho}_b=r\tilde{\rho}_{b0}(r),\ (r\tilde{\rho}_b)_t=-\lambda_b(r\tilde{\tilde{w}}_{b0})_{rr}(r),
\end{aligned}\right.
\end{equation*}
\begin{equation*}
\left\{\begin{aligned}
&(\partial_t^2-\lambda_b^2\partial_r^2)(r\tilde{\tilde{w}}_b)=G_b(t,r)\\
&t=0: r\tilde{\tilde{w}}_b=r\tilde{\tilde{w}}_{b0}(r),\ (r\tilde{\tilde{w}}_b)_t=-\lambda_br\tilde{\rho}_{b0}(r),
\end{aligned}\right.
\end{equation*}
then, we also have
\begin{align}
&\sum\limits_{\alpha=0,1}\|L_0^\alpha r\partial_r[(\tilde{\rho}_a+\tilde{w}_a)(\tilde{\rho}_b+\tilde{w}_b)]\|_{L^1([0,t]\times \mathbb{R}_+)}\notag\\
 \leq&C\Big(\|\tilde{\rho}_{a0}\|_{\dot{W}^{2,1}(\mathbb{R}^3)}+\|\tilde{w}_{a0}\|_{\dot{W}^{2,1}(\mathbb{R}^3)}+\sum\limits_{\alpha=0,1}\big(\|L_0^\alpha r\partial_rF_a\|_{L^1([0,t]\times\mathbb{R}_+)}+\|L_0^\alpha r\partial_r G_a\|_{L^1([0,t]\times\mathbb{R}_+)}\big)\Big)\notag\\
&\cdot\Big(\|\tilde{\rho}_{b0}\|_{\dot{W}^{2,1}(\mathbb{R}^3)}+\|\tilde{w}_{b0}\|_{\dot{W}^{2,1}(\mathbb{R}^3)}+\sum\limits_{\alpha=0,1}\big(\|L_0^\alpha F_b\|_{L^1([0,t]\times\mathbb{R}_+)}+\|L_0^\alpha G_b\|_{L^1([0,t]\times\mathbb{R}_+)}\big)\Big).
\end{align}

\end{thm}

\begin{rem}
Theorem \ref{thm2} and Theorem \ref{thm3} are natural generalizations of the following well-known bilinear estimate in one space dimension
\begin{equation*}
\left\{\begin{aligned}
&\phi_{it}+(\lambda_i \phi_i)_x=F_i(t,x)\\
&t=0: \phi_i=\phi_0^i(x),
\end{aligned}\right.
\end{equation*}
$i=1,2$. If $\lambda_1\neq\lambda_2$, then
\begin{align*}
\|\phi_1\phi_2\|_{L^1([0,t]\times \mathbb{R})}\leq C\Big(\|\phi_0^1\|_{L^1(\mathbb{R})}+\|F_1\|_{L^1([0,t]\times\mathbb{R})}\Big)\Big(\|\phi_0^2\|_{L^1(\mathbb{R})}+\|F_2\|_{L^1([0,t]\times\mathbb{R})}\Big).
\end{align*}
But the three dimensional estimate is much more difficult than the case of one space variable.
\end{rem}

It is worthy mention that the solution $u(t,\cdot)$ is not in the space $\dot{W}^{2,1}(R^3)$ for any time $t>0$, although $u_0\in\dot{W}^{2,1}(R^3)$, rather in some weighted Sobolev space
\begin{align*}
\sum\limits_{\alpha=0,1}\Big[\|L_0^\alpha\partial_t(ru(t,\cdot))\|_{L^1(\mathbb{R}_+)}+\|L_0^\alpha\partial_r(ru(t,\cdot))\|_{L^1(\mathbb{R}_+)}\Big]\leq C\varepsilon.
\end{align*}

Let us mention the existence theory of nonlinear wave equations given by Klainerman \cite{5}. He proved the global existence of classical solutions for 3D scalar quasilinear wave equation based on two basic assumptions, smallness of the initial data and null condition of the nonlinearities. This crucial work was also obtained independently by Christodoulou \cite{1} using a conformal mapping method. However, both of their work depend on the decay of corresponding linear wave equations. In this paper, we develop this to the hyperbolic systems which do not have any time decay.

Our long term project is trying to understand the formation of singularities and the global existence of weak solutions for systems of conservation laws in three dimensions. The multi-dimensional characteristics are so complicated that the tools which form the basis of a theory for hyperbolic conservation laws in a single space dimension do not extend to higher dimensions. Toward this goal, a key intermediate step is to understand solutions with radial symmetry exactly as that pointed out by Bressan et al in their book \cite{BA}. Unfortunately, radial symmetric systems do not admit global BV solutions generally (see \cite{BA, 8}), the theoretical analysis of several space dimensions still remains a grand challenge. Fortunately, we discover that the semilinear hyperbolic system allows global radial solutions with a bounded weighted BV norm. Thus to research the difficult hyperbolic systems in three dimensions, we consider the case with radial solutions as a first step.

Another motivation for doing this paper comes from the one dimensional quasilinear hyperbolic system
\begin{align}\label{1.9}
  u_t+A(u)u_x=0,
\end{align}
where A(u) is an n-th order square matrix. Let $\lambda_1(u),\cdots,\lambda_n(u)$ be the eigenvalues of $A(u)$, and $l_1(u),\cdots,l_n(u)$ denote its corresponding left eigenvectors. According to the important formula on the decomposition of waves which was given by John \cite{2}, if $w_i$ denotes the i-th component of $u_x$ in (\ref{1.9}),
\begin{align*}
 w_i=l_i(u)u_x,
\end{align*}
then we have
\begin{align}\label{1.10}
w_{it}+(\lambda_i(u)w_i)_x=\sum_{j\neq k}\gamma_{ijk}(u)w_jw_k.
\end{align}
Once we regard $\lambda_i(u)$ and $\gamma_{ijk}(u)$ as constants, (\ref{1.10}) is a semilinear system satisfying null condition.

For the completeness and for our future work, we also study the problem of quasilinear system in three dimensions. Consider the following general first-order quasilinear hyperbolic system
\begin{align}\label{2.8}
\frac{\partial u}{\partial t}+\sum_{i=1}^3A_i(u)\frac{\partial u}{\partial x^i}=0,
\end{align}
where  $u=(\rho_1,\cdots,\rho_l, v_1,\cdots,v_m)^T$ is the unknown function of $(t,x^1,x^2,x^3)$, $\rho_1,\cdots,\rho_l$ are scale unknown functions valued in $\mathbb{R}$, and $v_1,\cdots,v_m$ are vector unknown functions valued in $\mathbb{R}^3$. $A_1(u)$, $A_2(u)$, $A_3(u)$ are $(l+3m)\times(l+3m)$ matrixes with suitably smooth elements of $u$. We have
\begin{thm}\label{thm5}
If the general three dimensional quasilinear system (\ref{2.8}) is invariant under rotation of coordinates, and $\rho_k(t,x)=\rho_k(t,r)$, $v_p(t,x)=\frac{x}{r}w_p(t,r)$ are its radial solutions, then there exist $b_{kq}(\rho,w)$ and $c_{pj}(\rho,w)$ such that
\begin{equation}\left\{\begin{aligned}\label{2.12}
&\rho_{kt}+\sum\limits_{p=1}^mb_{kp}(\rho,w)w_{pr}+\sum\limits_{p=1}^m\frac{2}{r}b_{kp}(\rho,w)w_p=0\\
&w_{pt}+\sum\limits_{j=1}^lc_{pj}(\rho,w)\rho_{jr}=0.
\end{aligned}\right.\end{equation}
\end{thm}
Before ending the introduction, let us mention some important related works on classical solutions theory of quasilinear hyperbolic systems. In one spatial dimension case, the first result was given by John \cite{2}, he studied the formation of singularities when the system is genuinely nonlinear in the sense of Lax \cite{LPD}. Liu \cite{LTP} generalized John's result in a neighborhood of $u=0$, in which some linearly degenerate characteristic fields are allowed while other characteristics are genuinely nonlinear. H\"{o}rmander \cite{HL} reproved the result given in \cite{2} and obtained the sharp estimate for the life-span of $C^2$ solutions. By introducing the concept of weak linear degeneracy in \cite{6}, the global existence and the life-span of $C^1$ solutions to the general hyperbolic systems for small initial data were completed. More results on this problem have been established for linearly degenerate characteristics or weakly linearly degenerate characteristics with different smallness assumptions on the initial data, such as \cite{BA1, LZK2, ZY}. However, the the corresponding theory for general one-order multi-dimensional hyperbolic systems is still open, because of the complexity of characteristics. We also refer a part of outstanding works for the existence theory of special hyperbolic systems in multiple dimensions. For 3D nonlinear wave equations, John \cite{JF} and Alinhac \cite{AS} proved that the solution of Cauchy problem with sufficiently small initial data but disobeying null condition will blow up at a finite time. John and Klainerman \cite{JK} showed the almost global existence theory for nonlinear scalar wave equation. Sideris and Tu \cite{ST} developed Klainerman's generalized energy method in \cite{4} to obtain global existence under null condition. For compressible Euler equations, Sideris \cite{STC} considered the formation of singularities for a polytropic ideal fluid in 3D, where the initial data are constant outside a compact set. and a complete analysis of the effect of damping on the regularity and large time behavior of smooth solutions was done in \cite{STW}. By constructing special explicit solutions with spherical symmetry, Li and Wang \cite{LW} investigated the blowup phenomena for compressible Euler equations.

The remaining part of this paper is organized as follows. We devote section 2 to calculating the conditions which coefficient matrixes should satisfy, when the system (\ref{1.1}) is rotational invariant. And we give the proof of Theorem \ref{thm5}. We establish the null condition in section 3, and system (\ref{1.1}) is simplified under conditions [H1] and [H2]. In section 4, we derive the linear estimate which is needed. In section 5, we introduce a core bilinear estimate with non zero eigenvalues. Then we perform the bilinear estimate with zero eigenvalues in section 6.

\section{Rotational invariant system}
In this section, we study the general semilinear hyperbolic system (\ref{1.1}) which is invariant under rotation of coordinates, and we will give the proof of Lemma \ref{lem1} in the following.

In the case of rotational invariance, we shall treat scalar functions and vector functions separately. We set
\begin{align*}
  A_i=\begin{pmatrix}
           A^i & B^i \\
           C^i & D^i \\
        \end{pmatrix},
\end{align*}
where $A^i=(a_{jk}^i)$, $B^i=({\bf b}_{jq}^i)$, $C^i=({\bf c}_{pk}^i)$, $D^i=({\bf d}_{pq}^i)$, in which every element ${\bf b}_{jq}^i$ is a three dimensional row vector, ${\bf c}_{pk}^i$ is a three dimensional column vector, and ${\bf d}_{pq}^i$ is a $3\times3$ matrix, $i=1,2,3$; $j,k=1,\cdots,l$; $p,q=1,\cdots,m$.
To be specific,
\begin{align*}
   {\bf b}_{jq}^i=(b_{jq}^{i1},b_{jq}^{i2},b_{jq}^{i3}),\quad {\bf c}_{pk}^i=(c_{pk}^{i1},c_{pk}^{i2},c_{pk}^{i3})^T,\quad {\bf d}_{pq}^i=\begin{pmatrix}
                                                                                                                            d_{pq}^{i11} & d_{pq}^{i12} & d_{pq}^{i13} \\
                                                                                                                            d_{pq}^{i21} & d_{pq}^{i22} & d_{pq}^{i23} \\
                                                                                                                            d_{pq}^{i31} & d_{pq}^{i32} & d_{pq}^{i33} \\
                                                                                                                          \end{pmatrix}.
\end{align*}
Suppose that the quadratic term of (\ref{1.1}) has the following general formation:
\begin{equation*}
\begin{aligned}
Q_k^\rho(u,u)&=\sum_{i,j=1}^l\Gamma_{ij}^k\rho_i\rho_j+\sum_{1\leq j\leq l \atop 1\leq q\leq m}(\Upsilon_{jq}^{k1},\Upsilon_{jq}^{k2},\Upsilon_{jq}^{k3})\rho_j v_q+\sum_{s,q=1}^m
v_s^T \begin{pmatrix}
  \Omega_{sq}^{k11} & \Omega_{sq}^{k12} & \Omega_{sq}^{k13} \\
  \Omega_{sq}^{k21} & \Omega_{sq}^{k22} & \Omega_{sq}^{k23} \\
  \Omega_{sq}^{k31} & \Omega_{sq}^{k32} & \Omega_{sq}^{k33} \\
    \end{pmatrix}v_q,\\
Q_p^{v \alpha}(u,u)&=\sum_{i,j=1}^l \bar{\Gamma}_{ij}^{p\alpha}\rho_i\rho_j+\sum_{1\leq j\leq l \atop 1\leq q\leq m}(\bar{\Upsilon}_{jq}^{p\alpha1},\bar{\Upsilon}_{jq}^{p\alpha2},\bar{\Upsilon}_{jq}^{p\alpha3})\rho_j v_q\\
&\qquad\qquad\qquad\quad+\sum_{s,q=1}^m
v_s^T\begin{pmatrix}
  \bar{\Omega}_{sq}^{p\alpha11} & \bar{\Omega}_{sq}^{p\alpha12} & \bar{\Omega}_{sq}^{p\alpha13} \\
  \bar{\Omega}_{sq}^{p\alpha21} & \bar{\Omega}_{sq}^{p\alpha22} & \bar{\Omega}_{sq}^{p\alpha23} \\
  \bar{\Omega}_{sq}^{p\alpha31} & \bar{\Omega}_{sq}^{p\alpha32} & \bar{\Omega}_{sq}^{p\alpha33} \\
    \end{pmatrix}v_q,
    \end{aligned}\end{equation*}
where $\Gamma_{ij}^k$, $\Upsilon_{jq}^{k\beta}$, $\Omega_{sq}^{k\beta\gamma}$, $\bar{\Gamma}_{ij}^{p\beta}$, $\bar{\Upsilon}_{jq}^{p\alpha\beta}$, $\bar{\Omega}_{sq}^{p\alpha\beta\gamma}$, and $k=1,\cdots,l$; $p=1,\cdots,m$; $\alpha,\ \beta,\ \gamma=1,2,3.$

{\bf Proof of Lemma \ref{lem1}.}
Set $ {\bf O}=({\bf o}_1,{\bf o}_2,{\bf o}_3)^T=(o_{ij})\in \mathbb{R}^{3\times3}$ is any orthogonal matrix with $\det{\bf O}=1$, and let $y={\bf O}x$. Under rotation of coordinates, $\rho_k(t,x)$ and $v_p(t,x)$ become
\begin{align}\label{2.1}
  \rho_k(t,{\bf O}x),\ {\bf O}^Tv_p(t,{\bf O}x),
\end{align}
$k=1,\cdots,l$, $p=1,\cdots,m.$ For simplicity, we introduce two notations,
\begin{align*}
{\bf O}^T[v]&\triangleq({\bf O}^Tv_1,\cdots,{\bf O}^Tv_m)^T\\
{\bf d}({\bf \xi})&\triangleq {\bf d}^1{\bf \xi}_1+{\bf d}^2{\bf \xi}_2+{\bf d}^3{\bf \xi}_3,\ {\bf \xi}\in \mathbb{R}^3.
\end{align*}
Then take (\ref{2.1}) into system (\ref{1.1}), the left side becomes
\begin{align*}
\partial_t\begin{pmatrix}\rho(t,{\bf O}x)\\{\bf O}^T[v(t,{\bf O}x)] \end{pmatrix}+\sum\limits_{i=1}^3\begin{pmatrix} A^i & B^i\\C^i & D^i\end{pmatrix}\partial_{x^i}\begin{pmatrix}\rho(t,{\bf O}x)\\{\bf O}^T[v(t,{\bf O}x)] \end{pmatrix}.
\end{align*}
Thus we have
\begin{align}\label{2.3}
  \partial_t\begin{pmatrix}\rho(t,y)\\v(t,y) \end{pmatrix}&+\sum\limits_{i=1}^3\begin{pmatrix} \sum\limits_{j=1}^3A^jo_{ij} &  \sum\limits_{j=1}^3[B^j]{\bf O}^To_{ij}\\  \sum\limits_{j=1}^3{\bf O}[C^j]o_{ij} &  \sum\limits_{j=1}^3{\bf O}[ D^j]{\bf O}^To_{ij}\end{pmatrix}\partial_{y^i}\begin{pmatrix} \rho(t,y)\\v(t,y) \end{pmatrix}.
\end{align}
By assumption [H1], before and after coordinate rotation, the unknown functions should meet the same equation. So we obtain the following relationships,
$$\begin{pmatrix} a_{jk}^1 \\ a_{jk}^2 \\ a_{jk}^3 \end{pmatrix}={\bf O}\begin{pmatrix} a_{jk}^1 \\ a_{jk}^2 \\ a_{jk}^3 \end{pmatrix},\quad\begin{pmatrix} {\bf b}_{jq}^1 \\ {\bf b}_{jq}^2 \\ {\bf b}_{jq}^3 \end{pmatrix}={\bf O}\begin{pmatrix} {\bf b}_{jq}^1 \\ {\bf b}_{jq}^2 \\ {\bf b}_{jq}^3 \end{pmatrix}{\bf O}^T,$$
$$({\bf c}_{pk}^1,{\bf c}_{pk}^2,{\bf c}_{pk}^3)={\bf O}({\bf c}_{pk}^1,{\bf c}_{pk}^2,{\bf c}_{pk}^3){\bf O}^T,$$
$${\bf d}_{pq}^1={\bf O}{\bf d}_{pq}({\bf o}_1){\bf O}^T,\quad{\bf d}_{pq}^2={\bf O}{\bf d}_{pq}({\bf o}_2){\bf O}^T,\quad{\bf d}_{pq}^3={\bf O}{\bf d}_{pq}({\bf o}_3){\bf O}^T.$$
Since ${\bf O}$ is arbitrary, there are
\begin{align*}
\begin{pmatrix} a_{jk}^1 \\ a_{jk}^2 \\ a_{jk}^3 \end{pmatrix}=0,\quad\begin{pmatrix} {\bf b}_{jq}^1 \\ {\bf b}_{jq}^2 \\ {\bf b}_{jq}^3 \end{pmatrix}=\begin{pmatrix} b_{jq}^{11} &0&0\\ 0& b_{jq}^{11} &0 \\ 0&0& b_{jq}^{11} \end{pmatrix},\quad({\bf c}_{pk}^1,{\bf c}_{pk}^2,{\bf c}_{pk}^3)=\begin{pmatrix} c_{pk}^{11} &0&0\\ 0& c_{pk}^{11} &0 \\ 0&0& c_{pk}^{11} \end{pmatrix},\\
{\bf d}_{pq}^1=\begin{pmatrix}0&0&0\\ 0&0& d_{pq}^{111} \\ 0&-d_{pq}^{111} &0\end{pmatrix},\quad{\bf d}_{pq}^2=\begin{pmatrix} 0&0&-d_{pq}^{111}\\ 0&0&0\\ d_{pq}^{111}&0&0\end{pmatrix},\quad{\bf d}_{pq}^3=\begin{pmatrix} 0&d_{pq}^{111}&0\\-d_{pq}^{111}&0&0\\0&0&0\end{pmatrix},
\end{align*}
which gives $A_1$, $A_2$ and $A_3$ in Lemma \ref{lem1}.

The nonlinear part of (\ref{1.1}) can be treated in a similar way. We have algebraic conditions as follows,
  \begin{align*}
   (\Upsilon_{jq}^{k1},\Upsilon_{jq}^{k2},\Upsilon_{jq}^{k3})&=(\Upsilon_{jq}^{k1},\Upsilon_{jq}^{k2},\Upsilon_{jq}^{k3}){\bf O}^T,\\ \begin{pmatrix}
  \Omega_{sq}^{k11} & \Omega_{sq}^{k12} & \Omega_{sq}^{k13} \\
  \Omega_{sq}^{k21} & \Omega_{sq}^{k22} & \Omega_{sq}^{k23} \\
  \Omega_{sq}^{k31} & \Omega_{sq}^{k32} & \Omega_{sq}^{k33} \\
    \end{pmatrix}&={\bf O} \begin{pmatrix}
  \Omega_{sq}^{k11} & \Omega_{sq}^{k12} & \Omega_{sq}^{k13} \\
  \Omega_{sq}^{k21} & \Omega_{sq}^{k22} & \Omega_{sq}^{k23} \\
  \Omega_{sq}^{k31} & \Omega_{sq}^{k32} & \Omega_{sq}^{k33} \\
    \end{pmatrix}{\bf O}^T,
         \end{align*}
  \begin{gather*}
   \begin{pmatrix}
   \bar{\Gamma}_{ij}^{p1} \\
    \bar{\Gamma}_{ij}^{p2} \\
    \bar{\Gamma}_{ij}^{p3} \\
  \end{pmatrix}={\bf O}\begin{pmatrix}
    \bar{\Gamma}_{ij}^{p1} \\
    \bar{\Gamma}_{ij}^{p2} \\
    \bar{\Gamma}_{ij}^{p3} \\
  \end{pmatrix},\quad
   \begin{pmatrix}
   \bar{\Upsilon}_{jq}^{p11}&\bar{\Upsilon}_{jq}^{p12}&\bar{\Upsilon}_{jq}^{p13}\\                                                                                                                            \bar{\Upsilon}_{jq}^{p21}&\bar{\Upsilon}_{jq}^{p22}&\bar{\Upsilon}_{jq}^{p23}\\
   \bar{\Upsilon}_{jq}^{p31}&\bar{\Upsilon}_{jq}^{p32}&\bar{\Upsilon}_{jq}^{p33}\\
    \end{pmatrix}
   = {\bf O} \begin{pmatrix}
   \bar{\Upsilon}_{jq}^{p11}&\bar{\Upsilon}_{jq}^{p12}&\bar{\Upsilon}_{jq}^{p13}\\                                                                                                                            \bar{\Upsilon}_{jq}^{p21}&\bar{\Upsilon}_{jq}^{p22}&\bar{\Upsilon}_{jq}^{p23}\\
   \bar{\Upsilon}_{jq}^{p31}&\bar{\Upsilon}_{jq}^{p32}&\bar{\Upsilon}_{jq}^{p33}\\
    \end{pmatrix}{\bf O}^T,\\
   \begin{pmatrix}
   \bar{\Omega}_{sq}^{p\alpha11} & \bar{\Omega}_{sq}^{p\alpha12} & \bar{\Omega}_{sq}^{p\alpha13} \\
  \bar{\Omega}_{sq}^{p\alpha21} & \bar{\Omega}_{sq}^{p\alpha22} & \bar{\Omega}_{sq}^{p\alpha23} \\
  \bar{\Omega}_{sq}^{p\alpha31} & \bar{\Omega}_{sq}^{p\alpha32} & \bar{\Omega}_{sq}^{p\alpha33} \\
    \end{pmatrix}={\bf O}\sum_{\beta=1}^3\begin{pmatrix}
  \bar{\Omega}_{sq}^{p\beta11} & \bar{\Omega}_{sq}^{p\beta12} & \bar{\Omega}_{sq}^{p\beta13} \\
  \bar{\Omega}_{sq}^{p\beta21} & \bar{\Omega}_{sq}^{p\beta22} & \bar{\Omega}_{sq}^{p\beta23} \\
  \bar{\Omega}_{sq}^{p\beta31} & \bar{\Omega}_{sq}^{p\beta32} & \bar{\Omega}_{sq}^{p\beta33} \\
    \end{pmatrix}o_{\alpha\beta}{\bf O}^T.
 \end{gather*}
Then we can get
\begin{gather*}
   (\Upsilon_{jq}^{k1},\Upsilon_{jq}^{k2},\Upsilon_{jq}^{k3})=0,\quad
    \begin{pmatrix}
  \Omega_{sq}^{k11} & \Omega_{sq}^{k12} & \Omega_{sq}^{k13} \\
  \Omega_{sq}^{k21} & \Omega_{sq}^{k22} & \Omega_{sq}^{k23} \\
  \Omega_{sq}^{k31} & \Omega_{sq}^{k32} & \Omega_{sq}^{k33} \\
    \end{pmatrix}= \begin{pmatrix}
  \Omega_{sq}^{k11} & 0 & 0 \\
  0 & \Omega_{sq}^{k11} & 0 \\
  0 &0 & \Omega_{sq}^{k11} \\
    \end{pmatrix},\\
    \begin{pmatrix}
   \bar{\Gamma}_{ij}^{p1} \\
    \bar{\Gamma}_{ij}^{p2} \\
    \bar{\Gamma}_{ij}^{p3} \\
  \end{pmatrix}=0,\quad
  \begin{pmatrix}
   \bar{\Upsilon}_{jq}^{p11}&\bar{\Upsilon}_{jq}^{p12}&\bar{\Upsilon}_{jq}^{p13}\\                                                                                                                            \bar{\Upsilon}_{jq}^{p21}&\bar{\Upsilon}_{jq}^{p22}&\bar{\Upsilon}_{jq}^{p23}\\
   \bar{\Upsilon}_{jq}^{p31}&\bar{\Upsilon}_{jq}^{p32}&\bar{\Upsilon}_{jq}^{p33}\\
    \end{pmatrix}=\begin{pmatrix}
   \bar{\Upsilon}_{jq}^{p11}&0&0\\                                                                                                                            0&\bar{\Upsilon}_{jq}^{p11}&0\\
   0&0&\bar{\Upsilon}_{jq}^{p11}\\
    \end{pmatrix},\\
    \begin{pmatrix}
   \bar{\Omega}_{sq}^{p111} & \bar{\Omega}_{sq}^{p112} & \bar{\Omega}_{sq}^{p113} \\
  \bar{\Omega}_{sq}^{p121} & \bar{\Omega}_{sq}^{p122} & \bar{\Omega}_{sq}^{p123} \\
  \bar{\Omega}_{sq}^{p131} & \bar{\Omega}_{sq}^{p132} & \bar{\Omega}_{sq}^{p133} \\
    \end{pmatrix}=\begin{pmatrix}
   0 & 0 & 0 \\
  0& 0 & \bar{\Omega}_{sq}^{p111} \\
  0 & -\bar{\Omega}_{sq}^{p111} &0 \\
    \end{pmatrix},\\
    \begin{pmatrix}
   \bar{\Omega}_{sq}^{p211} & \bar{\Omega}_{sq}^{p212} & \bar{\Omega}_{sq}^{p213} \\
  \bar{\Omega}_{sq}^{p221} & \bar{\Omega}_{sq}^{p222} & \bar{\Omega}_{sq}^{p223} \\
  \bar{\Omega}_{sq}^{p231} & \bar{\Omega}_{sq}^{p232} & \bar{\Omega}_{sq}^{p233} \\
    \end{pmatrix}=\begin{pmatrix}
   0 & 0 & -\bar{\Omega}_{sq}^{p111} \\
  0& 0 &0 \\
  \bar{\Omega}_{sq}^{p111} & 0&0 \\
    \end{pmatrix},\\
    \begin{pmatrix}
   \bar{\Omega}_{sq}^{p311} & \bar{\Omega}_{sq}^{p312} & \bar{\Omega}_{sq}^{p313} \\
  \bar{\Omega}_{sq}^{p321} & \bar{\Omega}_{sq}^{p322} & \bar{\Omega}_{sq}^{p323} \\
  \bar{\Omega}_{sq}^{p331} & \bar{\Omega}_{sq}^{p332} & \bar{\Omega}_{sq}^{p333} \\
    \end{pmatrix}=\begin{pmatrix}
   0 & \bar{\Omega}_{sq}^{p111}& 0 \\
  -\bar{\Omega}_{sq}^{p111}& 0 &0 \\
  0 &0 &0 \\
    \end{pmatrix},
\end{gather*}
which gives that
\begin{align*}
Q_k^\rho(u,u)&=\sum_{i,j=1}^l\Gamma_{ij}^k\rho_i\rho_j+\sum_{s,q=1}^m \Omega_{sq}^{k11}v_s^T v_q,\\
Q_p^{v}(u,u)&=\sum_{1\leq j\leq l \atop 1\leq q\leq m} \bar{\Upsilon}_{jq}^{p11}\rho_j v_q+\sum_{s,q=1}^m \bar{\Omega}_{sq}^{p111}
v_s\times v_q.
\end{align*}
This has completed the proof of Lemma \ref{lem1}.\hfill$\Box$

In order to simplify the notations, set $b_{jq}\triangleq b_{jq}^{11}$, $c_{pk}\triangleq c_{pk}^{11}$, $d_{pq}\triangleq d_{pq}^{111}$. Take (\ref{1.5}) into system (\ref{1.1}), and using Lemma \ref{lem1}, we can easily get the following result by calculating directly.
\begin{thm}\label{thm4}
If system (\ref{1.1}) is rotational invariant and has radial solutions $\rho_k(t,x)=\rho_k(t,r)$, $v_p(t,x)=\frac{x}{r}w_p(t,r)$, then it can be reduced to
\begin{equation}\label{2.7}
\left\{\begin{aligned}
&\rho_t+Bw_r+\frac{2}{r}Bw=\sum\limits_{i,j=1}^l\Gamma_{ij}\rho_i\rho_j+\sum\limits_{s,q=1}^m\Omega_{sq}w_sw_q\\
&w_t+C\rho_r=\sum\limits_{1\leq j\leq l\atop 1\leq q\leq m}\Upsilon_{jq}\rho_jw_q,
\end{aligned}\right.
\end{equation}
where $\rho=(\rho_1,\cdots,\rho_l)^T$, $w=(w_1,\cdots,w_m)^T$, $B=(b_{jq})_{l\times m}$, $C=(c_{pk})_{m\times l}$, $\Gamma_{ij}=(\Gamma_{ij}^1,\cdots,\Gamma_{ij}^l)^T$, $\Omega_{sq}=(\Omega_{sq}^{111},\cdots,\Omega_{sq}^{l11})^T\triangleq(\Omega_{sq}^1,\cdots,\Omega_{sq}^l)^T$, $\Upsilon_{jq}=(\bar{\Upsilon}_{jq}^{l11},\cdots,\bar{\Upsilon}_{jq}^{m11})^T\triangleq(\Upsilon_{jq}^1,\cdots,\Upsilon_{jq}^m)^T$.
\end{thm}
Next, we will give a beginning of our future work and prove Theorem \ref{thm5}.

{\bf Proof of Theorem \ref{thm5}.} Similarly, let
\begin{align*}
A_i(u)=\begin{pmatrix}
           A^i(u) & B^i (u)\\
           C^i (u)& D^i(u) \\
        \end{pmatrix}
,\quad A^i=(a_{jk}^i(u)),\quad B^i=({\bf b}_{jq}^i(u)),\quad C^i=({\bf c}_{pk}^i(u)),D^i=({\bf d}_{pq}^i(u)),
\end{align*}
$i=1,2,3$; $j,k=1,\cdots,l$; $p,q=1,\cdots,m$, and
\begin{align*}
{\bf b}_{jq}^i(u)=(b_{jq}^{i1}(u),b_{jq}^{i2}(u),b_{jq}^{i3}(u)),\quad{\bf c}_{pk}^i(u)=(c_{pk}^{i1}(u),c_{pk}^{i2}(u),c_{pk}^{i3}(u))^T,
\end{align*}
\begin{align*}
{\bf d}_{pq}^i(u)=\begin{pmatrix}
                                       d_{pq}^{i11}(u) & d_{pq}^{i12}(u) & d_{pq}^{i13}(u) \\
                                       d_{pq}^{i21}(u) & d_{pq}^{i22}(u) & d_{pq}^{i23} (u)\\
                                       d_{pq}^{i31} (u)& d_{pq}^{i32}(u) & d_{pq}^{i33}(u) \\
                                       \end{pmatrix}.
\end{align*}

As the discussion in the semilinear case, if system (\ref{2.8}) is rotational invariant, then for any rotation matrix ${\bf O}$, $A_i(u)$ should satisfy
\begin{equation}\label{2.9}
\left.\begin{aligned}
\begin{pmatrix}
  a_{jk}^1(\rho,v) \\
  a_{jk}^2(\rho,v) \\
  a_{jk}^3 (\rho,v)\\
\end{pmatrix}={\bf O}\begin{pmatrix}
                           a_{jk}^1(\rho,{\bf O}^T[v]) \\
                           a_{jk}^2 (\rho,{\bf O}^T[v])\\
                           a_{jk}^3(\rho,{\bf O}^T[v]) \\
                         \end{pmatrix},\quad\begin{pmatrix}
                                                        {\bf b}_{jq}^1(\rho,v) \\
                                                        {\bf b}_{jq}^2(\rho,v) \\
                                                        {\bf b}_{jq}^3(\rho,v) \\
                                                      \end{pmatrix}={\bf O}\begin{pmatrix}
                                                                                 {\bf b}_{jq}^1 (\rho,{\bf O}^T[v])\\
                                                                                 {\bf b}_{jq}^2(\rho,{\bf O}^T[v]) \\
                                                                                 {\bf b}_{jq}^3(\rho,{\bf O}^T[v]) \\
                                                                               \end{pmatrix}{\bf O}^T,\\
({\bf c}_{pk}^1(\rho,v),{\bf c}_{pk}^2(\rho,v),{\bf c}_{pk}^3(\rho,v))={\bf O}({\bf c}_{pk}^1(\rho,{\bf O}^T[v]),{\bf c}_{pk}^2(\rho,{\bf O}^T[v]),{\bf c}_{pk}^3(\rho,{\bf O}^T[v])){\bf O}^T,\\
 {\bf d}_{pq}^1(\rho,v)={\bf O}({\bf d}_{pq}^1(\rho,{\bf O}^T[v])o_{11}+{\bf d}_{pq}^2(\rho,{\bf O}^T[v])o_{12}+{\bf d}_{pq}^3(\rho,{\bf O}^T[v])o_{13}){\bf O}^T,\\
{\bf d}_{pq}^2(\rho,v)={\bf O}({\bf d}_{pq}^1(\rho,{\bf O}^T[v])o_{21}+{\bf d}_{pq}^2(\rho,{\bf O}^T[v])o_{22}+{\bf d}_{pq}^3(\rho,{\bf O}^T[v])o_{23}){\bf O}^T,\\
{\bf d}_{pq}^3(\rho,v)={\bf O}({\bf d}_{pq}^1(\rho,{\bf O}^T[v])o_{31}+{\bf d}_{pq}^2(\rho,{\bf O}^T[v])o_{32}+{\bf d}_{pq}^3(\rho,{\bf O}^T[v])o_{33}){\bf O}^T.
\end{aligned}\right\}\end{equation}

Then we consider the radial solutions of the quasilinear system $\rho_k(t,x)=\rho_k(t,r)$, $v_p(t,x)=\frac{x}{r}w_p(t,r)$, and take them into (\ref{2.8}), we have
\begin{equation}\label{2.10}
\left\{\begin{aligned}
&\rho_{kt}+\sum\limits_{j=1}^l\sum\limits_{\alpha=1}^3\frac{x_\alpha}{r}a_{kj}^\alpha(\rho,v)\rho_{jr}+\sum\limits_{q=1}^m\sum\limits_{\alpha,\beta=1}^3\frac{x_\alpha x_\beta}{r^2}b_{kq}^{\alpha\beta}(\rho,v)w_{qr}\\
&\qquad-\sum\limits_{q=1}^m\sum\limits_{\alpha,\beta=1}^3\frac{x_\alpha x_\beta}{r^3}b_{kq}^{\alpha\beta}(\rho,v)w_q+\sum\limits_{q=1}^m\sum\limits_{\alpha=1}^3\frac{1}{r}b_{kq}^{\alpha\alpha}(\rho,v)w_q=0\\
&\frac{x}{r}w_{pt}+\sum\limits_{j=1}^l\sum\limits_{\alpha=1}^3\frac{x_\alpha}{r}{\bf c}_{pj}^\alpha(\rho,v)\rho_{jr}+\sum\limits_{q=1}^m\sum\limits_{\alpha,\beta=1}^3\frac{x_\alpha x_\beta}{r^2}{\bf d}_{pq}^{\alpha\beta}(\rho,v)w_{qr}\\
&\qquad-\sum\limits_{q=1}^m\sum\limits_{\alpha,\beta=1}^3\frac{x_\alpha x_\beta}{r^3}{\bf d}_{pq}^{\alpha\beta}(\rho,v)w_q
+\sum\limits_{q=1}^m\sum\limits_{\alpha=1}^3\frac{1}{r}{\bf d}_{pq}^{\alpha\alpha}(\rho,v)w_q=0,
\end{aligned}\right.\end{equation}
where $k=1,\cdots,l;p=1,\cdots,m,$ and ${\bf c}_{pj}^i=(c_{pj}^{i1},c_{pj}^{i1},c_{pj}^{i1})^T,\quad {\bf d}_{pq}^{\alpha\beta}=(d_{pq}^{\alpha1\beta},d_{pq}^{\alpha2\beta},d_{pq}^{\alpha3\beta})^T$.\\
Multiply $(\frac{x}{r})^T$ on both sides of the last equation, and take notice of (\ref{2.9}), we have
\begin{equation*}
\begin{aligned}
&\sum\limits_{\alpha=1}^3\frac{x_\alpha}{r}a_{kj}^\alpha(\rho,v)=a_{kj}^1(\rho,{\bf O}^T[v]),
\sum\limits_{\alpha,\beta=1}^3\frac{x_\alpha x_\beta}{r^2}b_{kq}^{\alpha\beta}(\rho,v)=b_{kq}^{11}(\rho,{\bf O}^T[v]),\\
&\sum\limits_{\alpha,\beta=1}^3\frac{x_\alpha x_\beta}{r^2}c_{pj}^{\alpha\beta}(\rho,v)=c_{pj}^{11}(\rho,{\bf O}^T[v]),
\sum\limits_{\alpha,\beta,\gamma=1}^3\frac{x_\alpha x_\beta x_\gamma}{r^3}d_{pq}^{\alpha\beta\gamma}(\rho,v)=d_{pq}^{111}(\rho,{\bf O}^T[v]).
\end{aligned}
\end{equation*}
So we can reduce (\ref{2.10}) to
\begin{equation}\label{2.11}
\left\{\begin{array}{ll}
\rho_{kt}+\sum\limits_{q=1}^mb_{kq}^{11}(\rho,{\bf O}^T[v])w_{qr}+\sum\limits_{q=1}^m\frac{2}{r}b_{kq}^{11}(\rho,{\bf O}^T[v])w_q=0\\
w_{pt}+\sum\limits_{j=1}^lc_{pj}^{11}(\rho,{\bf O}^T[v])\rho_{jr}=0.
\end{array}\right.
\end{equation}
Because ${\bf O}$ is arbitrary, there exists $b_{kq},c_{pj}$ such that $b_{kq}(\rho,w)=b_{kq}^{11}(\rho,{\bf O}^T[v]),c_{pj}(\rho,w)=c_{pj}^{11}(\rho,{\bf O}^T[v])$. \hfill$\Box$

\section{Null Condition}
We already know that $A_1$ and $Q(u,u)$ in (\ref{1.3}) have special formats in Lemma \ref{lem1}. In this section, we will study the nonlinear terms in (\ref{2.7}) under condition [H2]. If we rearrange the order of unknown functions in system (\ref{1.3}), let $V^{(\alpha)}=(v_1^\alpha,\cdots,v_m^\alpha)^T$, $\alpha=1,2,3$, then system (\ref{1.3}) equals to the following
\begin{align}\label{3.1}
\partial_t\begin{pmatrix}
    \rho \\
    V^{(1)}\\
    V^{(2)} \\
    V^{(3)} \\
  \end{pmatrix}+\begin{pmatrix}
                    0 & B & 0 & 0\\
                    C & 0 & 0 & 0\\
                    0 & 0 & D & 0\\
                    0 & 0 & 0 & -D\\
                  \end{pmatrix}\partial_1\begin{pmatrix}
                                 \rho \\
                                 V^{(1)}\\
                                 V^{(2)} \\
                                  V^{(3)} \\
                               \end{pmatrix}=\begin{pmatrix}
                                                     Q^\rho(u,u) \\
                                                     Q^{v1}(u,u) \\
                                                     Q^{v2}(u,u)\\
                                                     Q^{v3}(u,u)\\
                                                   \end{pmatrix},
\end{align}
where $B=(b_{jq})\in \mathbb{R}^{l\times m}$, $C=(c_{pk})\in \mathbb{R}^{m\times l}$, $D=(d_{pq})\in \mathbb{R}^{m\times m}$, $\bar{\Omega}_{sq}=(\bar{\Omega}_{sq}^{1111},\cdots,\bar{\Omega}_{sq}^{m111})^T$, and
\begin{align*}
 Q^\rho&=(Q^\rho_1,\cdots,Q^\rho_l),\ Q^{v\alpha}=(Q^{v\alpha}_1,\cdots,Q^{v\alpha}_m).
\end{align*}
Set
\begin{align*}
  M_1=\begin{pmatrix}
             0 & B\\
             C & 0 \\
             \end{pmatrix},\qquad M_2=\begin{pmatrix}
                                       D & 0 \\
                                       0 & -D \\
                                       \end{pmatrix}.
\end{align*}
Let $\Lambda=\text{diag}\{\lambda_1, \cdots,\lambda_{l+m}\}$ be the eigenvalue matrix of
$M_1$, $L=(l_{ij})_{(l+m)\times(l+m)}$ be its corresponding left  eigenvector matrix, and $R=(R_{ij})_{(l+m)\times (l+m)}=L^{-1}$.
Similarly, let $\tilde{L}M_2\tilde{R}=\tilde{\Lambda}$. These assumptions are reasonable because (\ref{1.3}) is a hyperbolic system.

For any fixed $\lambda_e$, $e\in\{1,\cdots,l+m\}$, consider the plane wave solutions to the linearized system of (\ref{3.1}),
\begin{align*}
\rho_k=\sum\limits_{a=1}^{l+m} R_{ka}Y_a(x-\lambda_e t),\ v_p^1=\sum\limits_{a=1}^{l+m} R_{l+p,a}Y_a(x-\lambda_e t),\\
v_p^2=\sum\limits_{b=1}^{2m}\tilde{R}_{pb}\tilde{Y}_b(x-\lambda_e t),\ v_p^3=\sum\limits_{b=1}^{2m}\tilde{R}_{m+p,b}\tilde{Y}_b(x-\lambda_e t),
\end{align*}
in which
$$Y_a(x-\lambda_e t)=\sum\limits_{j=1}^l l_{aj}\rho_j+\sum\limits_{q=1}^m l_{a,l+q}v_q^1,$$
$$\tilde{Y}_b(x-\lambda_e t)=\sum\limits_{q=1}^m \tilde{l}_{bq}v_q^2+\sum\limits_{q=1}^m l_{b,m+q}v_q^3.$$
Because of assumption [H2], we have
\begin{align}
\sum_{i,j=1}^l\Gamma_{ij}^kR_{ie}R_{je}+\sum_{s,q=1}^m\Omega_{sq}^kR_{l+s,e}R_{l+q,e}=0,\label{3.2}\\
\sum_{1\leq j\leq l\atop 1\leq q\leq m}\Upsilon_{jq}^{p}R_{je}R_{l+q,e}=0,\label{3.3}
\end{align}
where $k=1,\cdots,l$; $p=1,\cdots,m$.

 We remark that if the system is not strictly hyperbolic, without loss of generality, suppose that $\lambda_1=\lambda_2$, thus we furthermore obtain
\begin{align}\label{3.4}
\sum_{i,j=1}^l\Gamma_{ij}^kR_{i1}R_{j2}=\sum_{i,j=1}^l\Gamma_{ij}^kR_{i2}R_{j1}=0,
\end{align}
for all $k=1,\cdots,l$.

Next we pay attention to (\ref{2.7}). Rewrite $L=(L^1,L^2)$, $L^1\in \mathbb{R}^{(l+m)\times l}$, $L^2\in \mathbb{R}^{(l+m)\times m}$, multiply block $L^1$ on both sides of the first part of (\ref{2.7}) and multiply block $L^2$ on both sides of the second. Then we have
\begin{equation*}
\left\{\begin{aligned}
\tilde{\rho}_{et}+\lambda_e\tilde{w}_{er}+\frac{2}{r}\lambda_e\tilde{w}_a=&\sum\limits_{k=1}^l l_{ek}\sum\limits_{a,b=1}^{l+m}\Big[\sum\limits_{i,j=1}^l\Gamma_{ij}^kR_{ia}R_{jb}+\sum\limits_{s,q=1}^m\Gamma_{sq}^kR_{l+s,a}R_{l+q,b}\Big]\\
&\qquad\qquad\qquad\cdot(\tilde{\rho}_a+\tilde{w}_a)(\tilde{\rho}_b+\tilde{w}_b)\\
\tilde{w}_{et}+\lambda_{e}\tilde{\rho}_{er}=&\sum\limits_{p=1}^m l_{e,l+p}\sum\limits_{a,b=1}^{l+m}\sum\limits_{1\leq j\leq l \atop 1\leq q\leq m}\Upsilon_{jq}^pR_{ja}R_{l+q,b}(\tilde{\rho}_a+\tilde{w}_a)(\tilde{\rho}_b+\tilde{w}_b),
\end{aligned}\right.\end{equation*}
where
\begin{align}\label{3.5}
  \tilde{\rho}_a=\sum\limits_{j=1}^l l_{aj}\rho_j,\ \tilde{w}_a=\sum\limits_{q=1}^m l_{a,l+q}w_q.
\end{align}
Take notice of (\ref{3.2}) and (\ref{3.3}), we have
\begin{equation}\label{3.6}
\left\{\begin{aligned}
\tilde{\rho}_{et}+\lambda_e\tilde{w}_{er}+\frac{2}{r}\lambda_e\tilde{w}_e&=f_{1e}\\
\tilde{w}_{et}+\lambda_{e}\tilde{\rho}_{er}&=f_{2e},
\end{aligned}\right.\end{equation}
in which $f_{1e}$ and $f_{2e}$ are linear combinations of $(\tilde{\rho}_a+\tilde{w}_a)(\tilde{\rho}_b+\tilde{w}_b)$ and $e=1,\cdots,l+m$. It is crucial that $a\neq b$ by the algebraic condition (\ref{3.2}) and (\ref{3.3}). Furthermore, if the system has repeated eigenvalues $\lambda_1=\lambda_2$, then $b\neq 2$ $(a\neq 2)$ when $a=1$ $(b=1)$ because of (\ref{3.4}).

\section{Linear Estimate}
In the following contents of this article, we always use an uniform notation $f_e$ to represent for $f_{1e}$ and $f_{2e}$ appeared in (\ref{3.6}). We have
\begin{equation}\label{4.1}
\left\{\begin{aligned}
(\partial_t^2-\lambda_e^2\partial_r^2)(r\tilde{\rho}_e)&=r\partial_tf_e+r\partial_rf_e+f_e\\
(\partial_t^2-\lambda_e^2\partial_r^2)(r\tilde{\tilde{w}}_e)&=r\int\limits_0^r\partial_tf_ed\bar{r}+rf_e.
\end{aligned}\right.\end{equation}
But if $\lambda_e=0$, we should better replace (\ref{4.1}) by
\begin{equation}\label{4.2}
\left\{\begin{aligned}
\tilde{\rho}_{et}&=f_{1e}\\
\tilde{w}_{et}&=f_{2e}.
\end{aligned}\right.\end{equation}

Take notation $F$ as the nonlinear part of the first equation of (\ref{4.1}). If $\lambda_e\neq0$, we know that
\begin{align*}
\partial_t^2-\lambda_e^2\partial_r^2=(\partial_t\pm\lambda_e\partial_r)(\partial_t\mp\lambda_e\partial_r),
\end{align*}
and
\begin{align*}
[L_0,\Box_{\lambda_e}]=-2\Box_{\lambda_e},
\end{align*}
in which $\Box_{\lambda_e}$ denotes the wave operator with wave speed $\lambda_e$. Then
\begin{align*}
(\partial_t+\lambda_e\partial_r)(\partial_t-\lambda_e\partial_r)L_0^\alpha(r\tilde{\rho}_e)&=(L_0+2)^\alpha F,\\
(\partial_t-\lambda_e\partial_r)(\partial_t+\lambda_e\partial_r)L_0^\alpha(r\tilde{\rho}_e)&=(L_0+2)^\alpha F,
\end{align*}
where $L_0=t\partial_t+r\partial_r$, $\alpha=0,1$. The right side of the above equation can be written as $L_0^\alpha F$ briefly, this will have no essential difference in our analysis.

Let
\begin{align*}
  A=(\partial_t-\lambda_e\partial_r)L_0^\alpha(r\tilde{\rho}_e),\quad B=(\partial_t+\lambda_e\partial_r)L_0^\alpha(r\tilde{\rho}_e),
\end{align*}
then we have
\begin{align*}
 (\partial_t+\lambda_e\partial_r)|A|&=sgn(A)L_0^\alpha F,\\
 (\partial_t-\lambda_e\partial_r)|B|&=sgn(B)L_0^\alpha F.
\end{align*}
Thus
\begin{align*}
  \|A\|_{L^1(\mathbb{R}_+)}\leq \|A_0\|_{L^1(\mathbb{R}_+)}+\|L_0^\alpha F\|_{L^1([0,t]\times\mathbb{R}_+)}\leq C\|u_0\|_{\dot{W}^{2,1}(\mathbb{R}^3)}+\|L_0^\alpha F\|_{L^1([0,t]\times\mathbb{R}_+)},\\
   \|B\|_{L^1(\mathbb{R}_+)}\leq \|B_0\|_{L^1(\mathbb{R_+})}+\|L_0^\alpha F\|_{L^1([0,t]\times\mathbb{R}_+)}\leq C\|u_0\|_{\dot{W}^{2,1}(\mathbb{R}^3)}+\|L_0^\alpha F\|_{L^1([0,t]\times\mathbb{R}_+)}.
\end{align*}
That is
\begin{align*}
 \|\partial_tL_0^\alpha(r\tilde{\rho}_e)+\partial_rL_0^\alpha(r\tilde{\rho}_e)\|_{L^1(\mathbb{R}_+)}\leq C\|u_0\|_{\dot{W}^{2,1}(\mathbb{R}^3)}+\|L_0^\alpha F\|_{L^1([0,t]\times\mathbb{R}_+)}.
\end{align*}
And it is also true for $r\tilde{w}_e$, since
\begin{align*}
  r\tilde{w}_e=\partial_r(r\tilde{\tilde{w}}_e)-\tilde{\tilde{w}}_e.
\end{align*}
This implies
\begin{align}\label{4.3-0}
 \|\partial_tL_0^\alpha(ru)+\partial_rL_0^\alpha(ru)\|_{L^1(\mathbb{R}_+)}\leq C\|u_0\|_{\dot{W}^{2,1}(\mathbb{R}^3)}+\|L_0^\alpha F\|_{L^1([0,t]\times\mathbb{R}_+)}.
\end{align}

If $\lambda_e=0$, it is much easier to get (\ref{4.3-0}) because of (\ref{4.2}) and $[L_0,\partial_t]=-\partial_t$. We do not repeat that process here.

In order to obtain the global existence, by contraction mapping principle and global iterative method, we hope to have a uniform bound for $\tilde{\rho}_e$ and $\tilde{w}_e$ and their derivation with some weighted norm. The key to this is the bound for the quadratic terms.

\section{Bilinear Estimate with non zero eigenvalue}
Even if the nonlinear terms of (\ref{4.1}) have forms like $r\partial_tf_e$, $r\partial_rf_e$ and $f_e$, we only need to estimate $r\partial_rf_e$. Spirits of the proofs are similar, but the estimate of $r\partial_rf_e$ in some suitable weighted Sobolev space is much more complicated than others.

{\bf Proof of Theorem \ref{thm2}.}  We would estimate the bilinear estimate under the case of $\lambda_a=-\lambda_b$ separately. First, we pay attention to $|\lambda_a|\neq|\lambda_b|$.

For simplicity, according to Duhamel principle, we can only consider the following equations
\begin{equation*}
\left\{\begin{aligned}
&(\partial_t^2-\lambda_i^2\partial_r^2)(r\tilde{\rho}_i)=0\\
&t=0: \tilde{\rho}_i=0,\ \partial_t\tilde{\rho}_i=\varphi_i(r),
\end{aligned}\right.
\qquad
\left\{\begin{aligned}
&(\partial_t^2-\lambda_i^2\partial_r^2)(r\tilde{\tilde{w}}_i)=0\\
&t=0: \tilde{\tilde{w}}_i=0,\ \partial_t\tilde{\tilde{w}}_i=\psi_i(r),
\end{aligned}\right.\end{equation*}
where $i=a,b$. We would only study for $\alpha=0$ because of $[L_0,\Box]=-2\Box$. Thus, we need to perform the proof of the following
\begin{align}\label{4.3}
\|r\partial_r[(\tilde{\rho}_a+\tilde{w}_a)(\tilde{\rho}_b+\tilde{w}_b)]\|_{L^1([0,t]\times \mathbb{R}_+)} \leq& C\Big(\|s\varphi_{a}\|_{L^1(\mathbb{R}_+)}+\|(s\psi_{a})'\|_{L^1(\mathbb{R}_+)}\Big)\notag\\
&\cdot\Big(\|s^2(\varphi_{b})'\|_{L^1(\mathbb{R}_+)}+\|s(s\psi_{b})''\|_{L^1(\mathbb{R}_+)}\Big).
\end{align}
Since the wave speeds are positive and equals to $|\lambda_i|$, we can always assume $\lambda_i>0$ for notation convenience. If $\varphi_i$ and $\psi_i$ are spherical symmetric, we can express $\tilde{\rho}_i$ and $\tilde{\tilde{w}}_i$ in the form (see \cite{3})
\begin{align*}
  \tilde{\rho}_i(t,r)=\frac{1}{2\lambda_ir}\int_{|\lambda_it-r|}^{\lambda_it+r}s\varphi_i(s)ds,\  \tilde{\tilde{w}}_i(t,r)=\frac{1}{2\lambda_ir}\int_{|\lambda_it-r|}^{\lambda_it+r}s\psi_i(s)ds,
\end{align*}
then
\begin{align}\label{4.4}
 \tilde{\rho}_i+\tilde{w}_i=\frac{1}{r}\int_{|\lambda_it-r|}^{\lambda_it+r}sh_i(s)ds-\frac{1}{r^2}\int_{|\lambda_it-r|}^{\lambda_it+r}sm_i(s)ds,
\end{align}
where $sh_i(s)=\frac{1}{2\lambda_i}[s\varphi_i(s)+(s\psi_i(s))']$, $m_i(s)=\frac{1}{2\lambda_i}\psi_i(s)$.

Notice that
\begin{align}\label{4.5}
r\partial_r[(\tilde{\rho}_a+\tilde{w}_a)(\tilde{\rho}_b+\tilde{w}_b)](t,r)
=&\partial_r[r(\tilde{\rho}_a+\tilde{w}_a)]\cdot(\tilde{\rho}_b+\tilde{w}_b)-2(\tilde{\rho}_a+\tilde{w}_a)\cdot(\tilde{\rho}_b+\tilde{w}_b)\notag \\
&+(\tilde{\rho}_a+\tilde{w}_a)\cdot\partial_r[r(\tilde{\rho}_b+\tilde{w}_b)] \notag\\
\triangleq&I_A+I_B+I_C.
\end{align}

We first estimate $I_A$. Let $\varepsilon_0$ be small enough, we reduce matters to considering two cases:

{\bf Case 1:} $r\geq\varepsilon_0\tau.$

{\bf Case 2:} $r\leq\varepsilon_0\tau.$

Our strategy for dealing with $1/r$ is as follows. There is nonsingularity when $r$ is away from the origin, so we can treat is directly. If $r$ is near to the origin, we only need to contribute one derivative by Hadamard's formula.

{\bf Case 1:} By (\ref{4.4}) we have
\begin{align*}
\int_0^t\int_{\mathbb{R}_+}|I_A|drd\tau\leq&\int_0^t\int_{\mathbb{R}_+}\Big|(\lambda_a\tau+r)h_a(\lambda_a\tau+r)\frac{1}{r}\int_{|\lambda_b\tau-r|}^{\lambda_b\tau+r}sh_b(s)ds\Big|drd\tau\\
&+\int_0^t\int_{\mathbb{R}_+}\Big||\lambda_a\tau-r|h_a(|\lambda_a\tau-r|)\frac{1}{r}\int_{|\lambda_b\tau-r|}^{\lambda_b\tau+r}sh_b(s)ds\Big|drd\tau\\
&+\int_0^t\int_{\mathbb{R}_+}\Big|\frac{1}{r^2}\int_{|\lambda_a\tau-r|}^{\lambda_a\tau+r}(sm_a(s))'ds\int_{|\lambda_b\tau-r|}^{\lambda_b\tau+r}\bar{s}h_b(\bar{s})d\bar{s}\Big|drd\tau\\
&+\int_0^t\int_{\mathbb{R}_+}\Big|\frac{1}{r^3}\int_{|\lambda_a\tau-r|}^{\lambda_a\tau+r}sm_a(s)ds\int_{|\lambda_b\tau-r|}^{\lambda_b\tau+r}\bar{s}h_b(\bar{s})d\bar{s}\Big|drd\tau\\
&+\int_0^t\int_{\mathbb{R}_+}\Big|(\lambda_a\tau+r)h_a(\lambda_a\tau+r)\frac{1}{r^2}\int_{|\lambda_b\tau-r|}^{\lambda_b\tau+r}sm_b(s)ds\Big|drd\tau\\
&+\int_0^t\int_{\mathbb{R}_+}\Big||\lambda_a\tau-r|h_a(|\lambda_a\tau-r|)\frac{1}{r^2}\int_{|\lambda_b\tau-r|}^{\lambda_b\tau+r}sm_b(s)ds\Big|drd\tau\\
&+\int_0^t\int_{\mathbb{R}_+}\Big|\frac{1}{r^3}\int_{|\lambda_a\tau-r|}^{\lambda_a\tau+r}(sm_a(s))'ds\int_{|\lambda_b\tau-r|}^{\lambda_b\tau+r}\bar{s}m_b(\bar{s})d\bar{s}\Big|drd\tau\\
&+\int_0^t\int_{\mathbb{R}_+}\Big|\frac{1}{r^4}\int_{|\lambda_a\tau-r|}^{\lambda_a\tau+r}sm_a(s)ds\int_{|\lambda_b\tau-r|}^{\lambda_b\tau+r}\bar{s}m_b(\bar{s})d\bar{s}\Big|drd\tau\\
\triangleq&I_{A1}+I_{A2}+I_{A3}+I_{A4}+I_{A5}+I_{A6}+I_{A7}+I_{A8}.
\end{align*}

Let us considering $I_{A2}$ and $I_{A6}$ and then other terms are easy to follow, we estimate it in two subcases. For both cases we shall need to set $\xi=\lambda_a\tau-r$.

{\bf Subcase 1:} $\varepsilon_0\tau\leq r\leq \varepsilon_0^{-1}|\xi|$. Under this condition, we have $\frac{\varepsilon_0}{\lambda_a+\varepsilon_0}|\xi|\leq r\leq \varepsilon_0^{-1}|\xi|$. And from the integral expression of $I_{A6}$, we know that $0\leq s\leq (\frac{\lambda_b}{\varepsilon_0}+1)r$. Hence,
\begin{align*}
I_{A2}\leq& \|sh_b\|_{L^1(\mathbb{R}_+)}\int_0^t\int_{\varepsilon_0\tau\leq r\leq \varepsilon_0^{-1}|\xi|}|(\lambda_a\tau-r)h_a(|\lambda_a\tau-r|)|\cdot\frac{1}{r}drd\tau\\
\leq&C\|sh_b\|_{L^1(\mathbb{R}_+)}\int_{\mathbb{R}_+}|\xi h_a(\xi)|d\xi \int_{(\frac{\lambda_a}{\varepsilon_0}+1))^{-1}|\xi|}^{\varepsilon_0^{-1}|\xi|}\frac{1}{r}dr\\
\leq &C\|sh_b\|_{L^1(\mathbb{R}_+)}\|\xi h_a\|_{L^1(\mathbb{R}_+)},\\
I_{A6}\leq& \|sm_b\|_{L^\infty(\mathbb{R}_+)}\int_0^t\int_{\varepsilon_0\tau\leq r\leq \varepsilon_0^{-1}|\xi|}\big|(\lambda_a\tau-r)h_a(|\lambda_a\tau-r|)\frac{1}{r^2}\int_0^{(\frac{\lambda_b}{\varepsilon_0}+1)r}ds\big|drd\tau\\
\leq&C \|(sm_b)'\|_{L^1(\mathbb{R}_+)}\int_{\mathbb{R}_+}|\xi h_a(\xi)|d\xi\int_{(\frac{\lambda_a}{\varepsilon_0}+1)^{-1}|\xi|}^{\varepsilon_0^{-1}|\xi|}\frac{1}{r}dr\\
\leq &C\|(sm_b)'\|_{L^1(\mathbb{R}_+)}\|\xi h_a\|_{L^1(\mathbb{R}_+)}.
\end{align*}

{\bf Subcase 2:} $r\geq\varepsilon_0\tau$, and $r\geq \varepsilon_0^{-1}|\xi|$. In this case,
\begin{align*}
|\lambda_a-\lambda_b|r\leq\lambda_b|\lambda_a\tau -r|+\lambda_a|\lambda_b\tau-r|\leq\lambda_b\varepsilon_0r+\lambda_a|\lambda_b\tau-r|.
\end{align*}
Take notice of  $\lambda_a\neq\lambda_b$, $\lambda_a\lambda_b\neq0$, so $\varepsilon_0$ is sufficiently small and should be less than $\frac{|\lambda_a-\lambda_b|}{2\lambda_b}$, we have
\begin{align*}
\frac{|\lambda_a-\lambda_b|}{2\lambda_a}r\leq |\lambda_b\tau-r| \leq s\leq \lambda_b\tau+r\leq (\frac{\lambda_b}{\varepsilon_0}+1)r,
\end{align*}
which shows
\begin{align*}
\frac{\varepsilon_0}{\lambda_b+\varepsilon_0}s\leq r\leq\frac{2\lambda_a}{|\lambda_a-\lambda_b|}s.
\end{align*}
Hence,
\begin{align*}
I_{A2}\leq&\int_0^t\int_{r\geq\varepsilon_0\tau\atop r\geq \varepsilon_0^{-1}|\xi|}|(\lambda_a\tau-r)h_a(|\lambda_a\tau-r|)|\cdot|\frac{1}{r}\int_{\frac{|\lambda_a-\lambda_b|}{2\lambda_a}r}^{\frac{\lambda_b+\varepsilon_0}{\varepsilon_0}r}sh_b(s)ds|drd\tau\\
\leq&C\int_{\mathbb{R}_+}|\xi h_a(\xi)|d\xi\int_{\mathbb{R}_+}|sh_b(s)|ds\int_{\frac{\varepsilon_0}{\lambda_b+\varepsilon_0}s}^{\frac{2\lambda_a}{|\lambda_a-\lambda_b|}s}\frac{1}{r}dr\\
\leq &C\|\xi h_a\|_{L^1(\mathbb{R}_+)}\|sh_b\|_{L^1(\mathbb{R}_+)},
\end{align*}
\begin{align*}
I_{A6}\leq&\int_0^t\int_{r\geq\varepsilon_0\tau \atop r\geq \varepsilon_0^{-1}|\xi|}|(\lambda_a\tau-r)h_a(|\lambda_a\tau-r|)|\cdot|\frac{1}{r^2}\int_{\frac{|\lambda_a-\lambda_b|}{2\lambda_a}r}^{\frac{\lambda_b+\varepsilon_0}{\varepsilon_0}r}sm_b(s)ds|drd\tau\\
\leq&C\int_{\mathbb{R}_+}|\xi h_a(\xi)|d\xi\int_{\mathbb{R}_+}|sm_b(s)|ds\int_{\frac{\varepsilon_0}{\lambda_b+\varepsilon_0}s}^{\frac{2\lambda_a}{|\lambda_a-\lambda_b|}s}\frac{1}{r^2}dr\\
\leq&C\int_{\mathbb{R}_+}|\xi h_a(\xi)|d\xi\int_{\mathbb{R}_+}|sm_b(s)\frac{1}{s}|ds\\
\leq &C\|\xi h_a\|_{L^1(\mathbb{R}_+)}\|m_b\|_{L^1(\mathbb{R}_+)}.
\end{align*}

If we estimate $I_{A1}$ and $I_{A5}$ when $r\geq\varepsilon_0\tau$, we should set $\eta=\lambda_a\tau+r$. There naturally holds $r\leq \varepsilon_0^{-1}\eta$ as long as $\varepsilon_0$ is small enough. We omit the specific calculation here. And if we integral $I_{A3}$, $I_{A4}$ and $I_{A7}$, $I_{A8}$ by parts with respect to $r$, then it is same as $I_{A1}$, $I_{A2}$, $I_{A5}$ or $I_{A6}$. There is no essential difference with them.

{\bf Case 2:} Let $\xi_\theta=\lambda_a\tau+\theta r$ where $\theta\in[-1,1]$, since $r\leq\varepsilon_0\tau$, then
\begin{align*}
 \frac{1}{\lambda_a+\varepsilon_0}\xi_\theta\leq\tau\leq\frac{1}{\lambda_a-\varepsilon_0}\xi_\theta.
\end{align*}
We rewrite the integral of $|I_A|$ as
\begin{align*}
\int_0^t\int_{\mathbb{R}_+}|I_A|drd\tau\leq&\int_0^t\int_{\mathbb{R}_+}\Big|(\lambda_a\tau+r)h_a(\lambda_a\tau+r)\cdot(\tilde{\rho_b}+\partial_r\tilde{\tilde{w}}_b)\Big|drd\tau\\
&+\int_0^t\int_{\mathbb{R}_+}\Big|(\lambda_a\tau-r)h_a(\lambda_a\tau-r)\cdot(\tilde{\rho_b}+\partial_r\tilde{\tilde{w}}_b)\Big|drd\tau\\
&+\int_0^t\int_{\mathbb{R}_+}\Big|\frac{1}{r}\int_{|\lambda_a\tau-r|}^{\lambda_a\tau+r}(sm_a(s))'ds\cdot(\tilde{\rho_b}+\partial_r\tilde{\tilde{w}}_b)\Big|drd\tau\\
&+\int_0^t\int_{\mathbb{R}_+}\Big|\frac{1}{r^2}\int_{|\lambda_a\tau-r|}^{\lambda_a\tau+r}sm_a(s)ds\cdot(\tilde{\rho_b}+\partial_r\tilde{\tilde{w}}_b)\Big|drd\tau\\
\triangleq&I_{A9}+I_{A10}+I_{A11}+I_{A12}.
\end{align*}
Using $L^1-L^\infty$ estimate, which can be traced back to the 1960s, see Lemma 1 in reference \cite{7}, we have
\begin{align*}
 |\tilde{\rho_b}+\partial_r\tilde{\tilde{w}}_b|\leq C t^{-1}(\|r^2\varphi'_b\|_{L^1(\mathbb{R}_+)}+\|r^2\psi''_b\|_{L^1(\mathbb{R}_+)}).
\end{align*}
Thus,
\begin{align*}
I_{A9}\leq &C(\|r^2\varphi'_b\|_{L^1(\mathbb{R}_+)}+\|r^2\psi''_b\|_{L^1(\mathbb{R}_+)})\int_0^t\int_{r\leq\varepsilon_0\tau}|(\lambda_a\tau+r)h_a(\lambda_a\tau+r)\frac{1}{\tau}|drd\tau\\
\leq&C(\|r^2\varphi'_b\|_{L^1(\mathbb{R}_+)}+\|r^2\psi''_b\|_{L^1(\mathbb{R}_+)})\int_\mathbb{R}|\xi_1 h_a(\xi_1)|d\xi_1\int_{\frac{1}{\lambda_a+\varepsilon_0}\xi_1}^{\frac{1}{\lambda_a-\varepsilon_0}\xi_1}\frac{1}{\tau}d\tau\\
\leq&C\|\xi_1 h_a\|_{L^1(\mathbb{R}_+)}(\|r^2\varphi'_b\|_{L^1(\mathbb{R}_+)}+\|r^2\psi''_b\|_{L^1(\mathbb{R}_+)}),
\end{align*}
\begin{align*}
I_{A11}\leq &C(\|r^2\varphi'_b\|_{L^1(\mathbb{R}_+)}+\|r^2\psi''_b\|_{L^1(\mathbb{R}_+)})\\
&\cdot\int_0^t\int_{r\leq\varepsilon_0\tau}\frac{1}{r}\big|[(\lambda_a\tau+r)m_a(\lambda_a\tau+r)-(\lambda_a\tau-r)m_a(\lambda_a\tau-r)]\frac{1}{\tau}\big|drd\tau\\
\leq &C(\|r^2\varphi'_b\|_{L^1(\mathbb{R}_+)}+\|r^2\psi''_b\|_{L^1(\mathbb{R}_+)})\int_0^t\int_{r\leq\varepsilon_0\tau}\int_{-1}^1\Big|(sm_a)'\big|_{s=\lambda_a\tau+\theta r}\Big|d\theta \frac{1}{\tau}drd\tau\\
\leq&C(\|r^2\varphi'_b\|_{L^1(\mathbb{R}_+)}+\|r^2\psi''_b\|_{L^1(\mathbb{R}_+)})\int_{-1}^1d\theta\int_{\mathbb{R}_+}|(sm_a)'\big|_{s=\xi_\theta}|d\xi_\theta\int_{\frac{1}{\lambda_a+\varepsilon_0}\xi_\theta}^{\frac{1}{\lambda_a-\varepsilon_0}\xi_\theta}\frac{1}{\tau}d\tau\\
\leq&C\|(\xi_\theta m_a)'\|_{L^1(\mathbb{R}_+)}(\|r^2\varphi'_b\|_{L^1(\mathbb{R}_+)}+\|r^2\psi''_b\|_{L^1(\mathbb{R}_+)}).
\end{align*}

Obviously, the estimates of $I_{A10}$, $I_{A12}$ are also easy to follow because of $I_{A9}$, $I_{A11}$ and  the property of integration by parts. The method of treating $I_B$ and $I_C$ can be imitated as $I_A$. There is no essential difficulty, so we omit the detailed proof here.

Had not forgotten that, we have one last mission in this section, i.e. to prove the estimate when $\lambda_a=-\lambda_b\triangleq\lambda$. Without loss of generality, we assume that $\lambda>0$. According to the above proof, we discover that only subcase 2 in case 1 is highly dependent on $|\lambda_a|\neq|\lambda_b|$. Thus, we should separately consider the case for $r\geq\varepsilon_0\tau$ and $r\geq\varepsilon_0^{-1}|\xi|$ where $\xi=\lambda\tau-r$. Fortunately, we can get more information from the solutions of the following linear equations
\begin{equation*}
\left\{\begin{aligned}
&(\partial_t^2-\lambda^2\partial_r^2)(r\tilde{\rho}_a)=0\\
&t=0: r\tilde{\rho}_a=r\tilde{\rho}_{a0},\ \partial_tr\tilde{\rho}_a=-\lambda(r\tilde{\tilde{w}}_{a0})_{rr},
\end{aligned}\right.\qquad
\left\{\begin{aligned}
&(\partial_t^2-\lambda^2\partial_r^2)(r\tilde{\tilde{w}}_a)=0\\
&t=0: r\tilde{\tilde{w}}_a=r\tilde{\tilde{w}}_{a0},\ \partial_tr\tilde{\tilde{w}}_a=-\lambda r\tilde{\rho}_{a0},
\end{aligned}\right.
\end{equation*}
\begin{equation*}
\left\{\begin{aligned}
&(\partial_t^2-\lambda^2\partial_r^2)(r\tilde{\rho}_b)=0\\
&t=0: r\tilde{\rho}_b=r\tilde{\rho}_{b0},\ \partial_tr\tilde{\rho}_b=\lambda(r\tilde{\tilde{w}}_{b0})_{rr},
\end{aligned}\right.\qquad
\left\{\begin{aligned}
&(\partial_t^2-\lambda^2\partial_r^2)(r\tilde{\tilde{w}}_b)=0\\
&t=0: r\tilde{\tilde{w}}_b=r\tilde{\tilde{w}}_{b0},\ \partial_tr\tilde{\tilde{w}}_b=-\lambda r\tilde{\rho}_{b0},
\end{aligned}\right.\end{equation*}
then we can prove that
\begin{align}\label{4.8}
\|r\partial_r[(\tilde{\rho}_a+\tilde{w}_a)(\tilde{\rho}_b+\tilde{w}_b)]\|_{L^1([0,t]\times \mathbb{R}_+)} \leq& C\Big(\|r\tilde{\rho}_{a0}'\|_{L^1(\mathbb{R}_+)}+\|r\tilde{w}_{a0}'\|_{L^1(\mathbb{R}_+)}\Big)\notag\\
&\cdot\Big(\|\tilde{\rho}_{b0}\|_{\dot{W}^{2,1}(\mathbb{R}^3)}+\|\tilde{w}_{b0}\|_{\dot{W}^{2,1}(\mathbb{R}^3)}\Big).
\end{align}
Using D'Alembert formula, we have
\begin{align}
  r(\tilde{\rho}_a+\tilde{w}_a)=&(r-\lambda t)\tilde{\rho}_{a0}(r-\lambda t)+\tilde{\tilde{w}}_{a0}(r-\lambda t)+(r-\lambda t)\tilde{w}_{a0}(r-\lambda t)\notag \\
  &-\frac{1}{2r}[(r-\lambda t)\tilde{\tilde{w}}_{a0}(r-\lambda t)+(r+\lambda t)\tilde{\tilde{w}}_{a0}(r+\lambda t)]+\frac{1}{2r}\int_{|\lambda t-r|}^{\lambda t+r}s\tilde{\rho}_{a0}(s)ds,
\end{align}
\begin{align}\label{34}
   r(\tilde{\rho}_b+\tilde{w}_b)=&(r+\lambda t)\tilde{\rho}_{b0}(r+\lambda t)+\tilde{\tilde{w}}_{b0}(r+\lambda t)+(r+\lambda t)\tilde{w}_{b0}(r+\lambda t)\notag\\
  &-\frac{1}{2r}[(r+\lambda t)\tilde{\tilde{w}}_{b0}(r+\lambda t)+(r-\lambda t)\tilde{\tilde{w}}_{b0}(r-\lambda t)]+\frac{1}{2r}\int_{|\lambda t-r|}^{\lambda t+r}s\tilde{\rho}_{b0}(s)ds.
\end{align}
Then
\begin{align}\label{35}
\partial_r[r(\tilde{\rho}_a+\tilde{w}_a)]=&(r-\lambda t)\tilde{\rho}'_{a0}(r-\lambda t)+\tilde{\rho}_{a0}(r-\lambda t)+2\tilde{w}_{a0}(r-\lambda t)+(r-\lambda t)\tilde{w}'_{a0}(r-\lambda t)\notag\\
&-\frac{1}{2r}[(r-\lambda t)\tilde{w}_{a0}(r-\lambda t)+\tilde{\tilde{w}}_{a0}(r-\lambda t)+(r+\lambda t)\tilde{w}_{a0}(r+\lambda t)+\tilde{\tilde{w}}_{a0}(r+\lambda t)\notag\\
&\qquad+(r-\lambda t)\tilde{\rho}_{a0}(r-\lambda t)-(r+\lambda t)\tilde{\rho}_{a0}(r+\lambda t)]\notag\\
&+\frac{1}{2r^2}[(r-\lambda t)\tilde{\tilde{w}}_{a0}(r-\lambda t)+(r+\lambda t)\tilde{\tilde{w}}_{a0}(r+\lambda t)-\int_{|\lambda t-r|}^{\lambda t+r}s\tilde{\rho}_{a0}(s)ds].
\end{align}
We just show the proof of following three terms of $I_A=\partial_r[r(\tilde{\rho}_a+\tilde{w}_a)]\cdot(\tilde{\rho}_b+\tilde{w}_b)$, and the others can be estimated similarly. Since $r\geq\varepsilon_0\tau$ and $r\geq \varepsilon_0^{-1}|\xi|$, we have $r\leq\eta=r+\lambda\tau\leq(1+\lambda\varepsilon_0^{-1})r$. Inserting (\ref{34}) and (\ref{35}) into $I_A$, each term can be followed as one of the following
\begin{align*}
  &\int_0^t\int_{r\geq\varepsilon_0\tau\atop r\geq \varepsilon_0^{-1}|\xi|}\big|\tilde{\rho}_{a0}(|r-\lambda\tau|)\cdot\tilde{\tilde{w}}_{b0}(|r-\lambda\tau|)\frac{|r-\lambda\tau|}{r^2}\big|drd\tau\\
  \leq&C\int_{\mathbb{R}_+}\int_{r\geq \varepsilon_0^{-1}|\xi|}\big|\tilde{\rho}_{a0}(|\xi|)\tilde{\tilde{w}}_{b0}(|\xi|)\frac{\xi}{r^2}\big|drd\xi\\
 \leq&C\|\tilde{\tilde{w}}_{b0}\|_{L^\infty(\mathbb{R}_+)}\int_{\mathbb{R}_+}\big|\tilde{\rho}_{a0}(|\xi|)\frac{\xi}{\xi}\big|drd\xi\\
  \leq& C\|\tilde{\rho}_{a0}\|_{L^1(\mathbb{R}_+)}\|\tilde{w}_{b0}\|_{L^1(\mathbb{R}_+)},
\end{align*}
\begin{align*}
  &\int_0^t\int_{r\geq\varepsilon_0\tau\atop r\geq \varepsilon_0^{-1}|\xi|}\big|\frac{1}{r}\tilde{\tilde{w}}_{a0}(|r-\lambda\tau|)\cdot\frac{1}{r}\tilde{\tilde{w}}_{b0}(r+\lambda\tau)\big|drd\tau\\
  \leq&C\Big(\int_0^t\int_{r\geq\varepsilon_0\tau\atop r\geq \varepsilon_0^{-1}|\xi|}\big|\frac{1}{r^2}\tilde{\tilde{w}}_{a0}^2(|r-\lambda\tau|)\big|drd\tau\Big)^{\frac{1}{2}}\Big(\int_0^t\int_{r\geq(1+\frac{\lambda}{\varepsilon_0})^{-1}\eta}\big|\frac{1}{r^2}\tilde{\tilde{w}}_{b0}^2(r+\lambda\tau)\big|drd\tau\Big)^{\frac{1}{2}}\\
  \leq&C\Big(\int_{\mathbb{R}_+}\Big|\frac{(\int_0^{|\xi|}\tilde{w}_{a0}(s)ds)^2}{\xi}\Big|d\xi\Big)^{\frac{1}{2}}\Big(\int_{\mathbb{R}_+}\Big|\frac{(\int_0^\eta\tilde{w}_{b0}(s)ds)^2}{\eta}\Big|d\eta\Big)^{\frac{1}{2}}\\
  \leq&C\Big(\int_{\mathbb{R}_+}\int_0^{|\xi|}\tilde{w}_{a0}^2(s)dsd\xi\Big)^{\frac{1}{2}}\Big(\int_{\mathbb{R}_+}\int_0^\eta\tilde{w}_{b0}^2(s)dsd\eta\Big)^{\frac{1}{2}}\\
  \leq&C\Big(\int_{\mathbb{R}_+}|\xi|\tilde{w}_{a0}^2(\eta)d\xi\Big)^{\frac{1}{2}}\Big(\int_{\mathbb{R}_+}\eta\tilde{w}_{b0}^2(\eta)dsd\eta\Big)^{\frac{1}{2}}\\
  \leq&C\|(\xi\tilde{w}_{a0})'\|_{L^1(\mathbb{R}_+)}\|(\eta\tilde{w}_{b0})'\|_{L^1(\mathbb{R}_+)},
\end{align*}
\begin{align*}
&\int_0^t\int_{r\geq\varepsilon_0\tau\atop r\geq \varepsilon_0^{-1}|\xi|}\big|\frac{1}{r^2}\int_{|\lambda t-r|}^{\lambda t+r}\sigma\tilde{\rho}_{a0}(\sigma)d\sigma\frac{1}{r^2}\int_{|\lambda t-r|}^{\lambda t+r}s\tilde{\rho}_{b0}(s)ds\big|drd\tau\\
\leq&C\int_{\mathbb{R}_+}|\sigma\tilde{\rho}_{a0}(\sigma)|\int_{\mathbb{R}_+}|s\tilde{\rho}_{b0}(s)|\int_{r\geq(1+\frac{\lambda}{\varepsilon_0})^{-1}s\atop r\geq(1+\frac{\lambda}{\varepsilon_0})^{-1}\sigma}\frac{1}{r^4}\int_{\tau\leq\varepsilon_0^{-1}r}d\tau drdsd\sigma\\
\leq&C\|\tilde{\rho}_{a0}\|_{L^1(\mathbb{R}_+)}\|\tilde{\rho}_{b0}\|_{L^1(\mathbb{R}_+)}.
\end{align*}

We deal $I_B$ and $I_C$ with the same method like $I_A$ when $\lambda_a=-\lambda_b=\lambda$. By now the proof of Theorem {\ref{thm2}} is finished, since $\|sh_b\|_{L^1(\mathbb{R}_+)}\leq \|s^2h'_b\|_{L^1(\mathbb{R}_+)}$. \hfill$\Box$

\section{Bilinear Estimate with zero eigenvalue}
In this section, we consider situation with zero eigenvalue.

{\bf Proof of Theorem \ref{thm3}.} By Duhamel principle and along the discussion in section 5, we need to estimate
\begin{align}\label{5.0}
&\|r\partial_r[(\tilde{\rho}_a+\tilde{w}_a)(\tilde{\rho}_b+\tilde{w}_b)]\|_{L^1([0,t]\times \mathbb{R}_+)}\notag\\
 \leq & C\Big(\|r\tilde{\rho}'_{a0}\|_{L^1(\mathbb{R}_+)}+\|r\tilde{w}'_{a0}\|_{L^1(\mathbb{R}_+)}\Big)\cdot\Big(\|s^2\varphi'_{b}\|_{L^1(\mathbb{R}_+)}+\|s(s\psi_{b})''\|_{L^1(\mathbb{R}_+)}\Big).
\end{align}

From (\ref{5.0}), we know that $\lambda_a$ and $\lambda_b$ play different roles. This makes us think of dealing with $\lambda_a=0$ and $\lambda_b=0$ differently. Certainly, it is impossible that both of them equal to zero because of the null condition (\ref{3.2}) and (\ref{3.3}).

If $\lambda_a=0$, $\lambda_b>0$, then $\tilde{\rho}_a$, $\tilde{w}_a$ satisfy ordinary differential equations,
\begin{equation}\label{5.1}
\left\{\begin{aligned}
\tilde{\rho}_{at}&=f_{1a}\\
\tilde{w}_{at}&=f_{2a}.
\end{aligned}\right.\end{equation}
Using expression (\ref{4.4}), we have
\begin{align*}
  \|I_A\|_{L^1([0,t]\times \mathbb{R}_+)}\leq&\int_0^t\int_{\mathbb{R}_+}\big|r(\tilde{\rho}'_{a0}+\tilde{w}'_{a0})\frac{1}{r}\int_{|\lambda_b\tau-r|}^{\lambda_b\tau+r}sh_b(s)ds\big|drd\tau\\
  &+\int_0^t\int_{\mathbb{R}_+}\big|r(\tilde{\rho}'_{a0}+\tilde{w}'_{a0})\frac{1}{r^2}\int_{|\lambda_b\tau-r|}^{\lambda_b\tau+r}sm_b(s)ds\big|drd\tau\\
  \triangleq & I_{Ah}+I_{Am}.
\end{align*}
They can be estimated as
\begin{align*}
  I_{Ah}\leq&\int_0^t\int_{\mathbb{R}_+}\big|r(\tilde{\rho}'_{a0}+\tilde{w}'_{a0})\frac{1}{r}[H_b(\lambda_b\tau+r)-H_b(\lambda_b\tau-r)]\big|drd\tau\\
\leq&\int_0^t\int_{\mathbb{R}_+}\big|r(\tilde{\rho}'_{a0}+\tilde{w}'_{a0})\int_{-1}^1(\lambda_b\tau+\theta r)h_b(\lambda_b\tau+\theta r)d\theta\big|drd\tau\\
\leq &C \int_{\mathbb{R}_+}|r(\tilde{\rho}'_{a0}+\tilde{w}'_{a0})|dr\int_{\mathbb{R}_+}|\xi h_b(\xi)|dr,
\end{align*}
\begin{align*}
I_{Am}\leq&\int_{\mathbb{R}_+}|r(\tilde{\rho}'_{a0}+\tilde{w}'_{a0})|\cdot\int_0^t\big|\frac{1}{r^2}\int_{|\lambda_b\tau-r|}^{\lambda_b\tau+r}sm_b(s)dsd\tau\big|dr\\
\leq&ess\sup_{r\in \mathbb{R}_+} \Big|\int_0^t\frac{1}{r^2}\int_{|\lambda_b\tau-r|}^{\lambda_b\tau+r}sm_b(s)dsd\tau\Big|\cdot \int_{\mathbb{R}_+}|r(\tilde{\rho}'_{a0}+\tilde{w}'_{a0})|dr\\
\leq&C\|r(\tilde{\rho}'_{a0}+\tilde{w}'_{a0})\|_{L^1(\mathbb{R}_+)}\int_{\mathbb{R}_+}\int_0^t\Big|\frac{1}{r^3}\int_{|\lambda_b\tau-r|}^{\lambda_b\tau+r}sm_b(s)ds\Big|d\tau dr\\
 &+\|r(\tilde{\rho}'_{a0}+\tilde{w}'_{a0})\|_{L^1(\mathbb{R}_+)} \int_{\mathbb{R}_+}\int_0^t\Big|\frac{1}{r^2}\int_{|\lambda_b\tau-r|}^{\lambda_b\tau+r}(sm_b(s))'ds\Big|d\tau dr\\
 \leq & C\|r(\tilde{\rho}'_{a0}+\tilde{w}'_{a0})\|_{L^1(\mathbb{R}_+)} \int_0^t \int_{\mathbb{R}_+}\Big|\frac{1}{r^2}\int_{|\lambda_b\tau-r|}^{\lambda_b\tau+r}(sm_b(s))'ds \Big|drd\tau.
\end{align*}
Set
\begin{align*}
I_{Am1}&=\int_0^t\int_{\mathbb{R}_+}\Big|\frac{1}{r^2}\int_{|\lambda_b\tau-r|}^{\lambda_b\tau+r}(sm_b(s))'ds\Big|dr d\tau\\
&=\int_0^t\int_{\mathbb{R}_+}\frac{1}{r^2}|(\lambda_b\tau+r )m_b(\lambda_b\tau+r )-|\lambda_b\tau-r|m_b(|\lambda_b\tau-r|)|dr d\tau,
\end{align*}
we estimate it by considering two cases:

{\bf Case 1:} $r\geq\varepsilon_0\tau$. In this case, we have $\xi=\lambda_b\tau+r\leq(\frac{\lambda_b}{\varepsilon_0}+1)r$, $\eta=\lambda_b\tau-r\leq\lambda_b\tau+r\leq(\frac{\lambda_b}{\varepsilon_0}+1)r$, then
\begin{align*}
\int_0^t\int_{r\geq\varepsilon_0\tau}\frac{1}{r^2}|(\lambda_b\tau+r )m_b(\lambda_b\tau+r )|dr d\tau\leq& C\int_{\mathbb{R}_+}\int_{r\geq(\frac{\lambda_b}{\varepsilon_0}+1)^{-1}\xi}\frac{1}{r^2}|\xi m_b(\xi )|dr d\xi\\
\leq&C\int_{\mathbb{R}_+}|m_b(\xi )|d\xi,
\end{align*}
\begin{align*}
\int_0^t\int_{r\geq\varepsilon_0\tau}\frac{1}{r^2}|\lambda_b\tau-r||m_b(|\lambda_b\tau-r||)|dr d\tau\leq& C\int_{\mathbb{R}_+}\int_{r\geq(\frac{\lambda_b}{\varepsilon_0}+1)^{-1}\eta}\frac{1}{r^2}|\eta m_b(\eta )|dr d\eta\\
\leq&C\int_{\mathbb{R}_+}|m_b(\eta )|d\eta.
\end{align*}

{\bf Case 2:} $r\leq\varepsilon_0\tau$. Then $\xi_\theta=\lambda_b\tau+\theta r\geq \lambda_b\tau-r\geq(\frac{\lambda_b}{\varepsilon_0}-1)r$, and
\begin{align*}
I_{Am1}&\leq\int_0^t\int_{r\leq\varepsilon_0\tau}\Big|\int_{-1}^1(sm_b(s))''|_{s=\lambda_b\tau+\theta r }d\theta \Big|dr d\tau\\
&\leq C\int_{-1}^1\int_{\mathbb{R}_+}\int_{r\leq(\frac{\lambda_b}{\varepsilon_0}-1)^{-1}\xi_\theta}|(\xi_\theta m_b(\xi_\theta))''|dr d\xi_\theta d\theta\leq C\int_{\mathbb{R}_+}|\xi_\theta(\xi_\theta m_b(\xi_\theta))''|d\xi_\theta.
\end{align*}

Since
\begin{align*}
  \|I_B\|_{L^1([0,t]\times \mathbb{R}_+)}\leq&\int_0^t\int_{\mathbb{R}_+}\big|(\tilde{\rho}_{a0}+\tilde{w}_{a0})\frac{1}{r}\int_{|\lambda_b\tau-r|}^{\lambda_b\tau+r}sh_b(s)ds\big|drd\tau\\
  &+\int_0^t\int_{\mathbb{R}_+}\big|(\tilde{\rho}_{a0}+\tilde{w}_{a0})\frac{1}{r^2}\int_{|\lambda_b\tau-r|}^{\lambda_b\tau+r}sm_b(s)ds\big|drd\tau,
  \end{align*}
then the estimate of $\|I_B\|_{L^1([0,t]\times \mathbb{R}_+)}$ is like $ \|I_A\|_{L^1([0,t]\times \mathbb{R}_+)}$ when $\lambda_a=0$, $\lambda_b\neq 0$.

 We calculate $\partial_r(r(\tilde{\rho}_{b}+\tilde{w}_{b}))$,
 \begin{align*}
  \partial_r(r(\tilde{\rho}_{b}+\tilde{w}_{b}))=&[(\lambda_b\tau+r)h_b(\lambda_b\tau+r)-(\lambda_b\tau-r)h_b(\lambda_b\tau-r)]\\
  &-\frac{1}{r}[(\lambda_b\tau+r)m_b(\lambda_b\tau+r)-(\lambda_b\tau-r)m_b(\lambda_b\tau-r)]+\frac{1}{r^2}\int_{|\lambda_b\tau-r|}^{\lambda_b\tau+r}sm_b(s)ds,
  \end{align*}
 thus we can estimate $\|I_C\|_{L^1([0,t]\times \mathbb{R}_+)}$ by the method what we have used as above.

If $\lambda_a>0$, $\lambda_b=0$, then $\tilde{\rho}_b$, $\tilde{w}_b$ satisfy ordinary differential equations, and $\tilde{\rho}_b(t,r)=\tilde{\rho}_{b0}(r)$, $\tilde{w}_b(t,r)=\tilde{w}_{b0}(r)$. From the previous discussion, we only need to pay attention to the case of $r\leq \varepsilon_0\tau$. We know
\begin{align}\label{5.2}
r\partial_r[(\tilde{\rho}_a+\tilde{w}_a)(\tilde{\rho}_b+\tilde{w}_b)]
=&\partial_r[r(\tilde{\rho}_a+\tilde{w}_a)]\cdot(\tilde{\rho}_b+\tilde{w}_b)+(\tilde{\rho}_a+\tilde{w}_a)\cdot r\partial_r(\tilde{\rho}_b+\tilde{w}_b)\notag\\
&-(\tilde{\rho}_a+\tilde{w}_a)(\tilde{\rho}_b+\tilde{w}_b),
\end{align}
and expression
\begin{align}\label{4.7}
  \tilde{\rho}_b+\tilde{w}_b=&\frac{1}{r}\int_{|\lambda_b t-r|}^{\lambda_bt+r}sg_b(s)ds+\partial_r\int_{-1}^1sm_b(s)\big|_{s=\lambda_bt+\theta r}d\theta \notag\\
  =&\int_{-1}^1sg_b(s)\big|_{s=\lambda_b t+\theta r}d\theta+\int_{-1}^1(sm_b)'\big|_{s=\lambda_b t+\theta r}\theta d\theta,
\end{align}
where $g_b=-\frac{1}{2\lambda}\varphi_b^{(2)}$, will help us to discuss.

We obtain
\begin{align*}
 &\int_0^t\int_{r\leq\varepsilon_0\tau}\Big|(\lambda_a\tau\pm r)h_a(\lambda_a\tau\pm r)\cdot(\tilde{\rho}_{b0}(r)+\tilde{w}_{b0}(r))\Big|drd\tau\\
 \leq &C\int_{\mathbb{R}_+}\int_{\mathbb{R}_+}|\xi_\pm h_a(\xi_\pm)\cdot(\tilde{\rho}_{b0}(r)+\tilde{w}_{b0}(r))|drd\xi_\pm\\
\leq &C\|\xi_\pm h_a(\xi_\pm)\|_{L^1(\mathbb{R}_+)}\|\tilde{\rho}_{b0}+\tilde{w}_{b0}\|_{L^1(\mathbb{R}_+)},
\end{align*}
\begin{align*}
 &\int_0^t\int_{r\leq\varepsilon_0\tau}\Big|\int_{-1}^1(sm_a)'\big|_{s=\lambda_a\tau+\theta r}\theta d\theta\cdot(\tilde{\rho}_{b0}(r)+\tilde{w}_{b0}(r))\Big|drd\tau\\
\leq &C\int_{-1}^1\int_{\mathbb{R}_+}\int_{\mathbb{R}_+}|(\xi_\theta m_a(\xi_\theta))'\cdot(\tilde{\rho}_{b0}(r)+\tilde{w}_{b0}(r))|drd\xi_\theta d\theta\\
\leq& C\|(\xi_\theta m_a)'\|_{L^1(\mathbb{R}_+)}\|\tilde{\rho}_{b0}+\tilde{w}_{b0}\|_{L^1(\mathbb{R}_+)},
\end{align*}
where $\xi_\pm=\lambda_a\tau\pm r$, $\xi_\theta=\lambda_a\tau+\theta r$. Thus
%\begin{align*}
%  \tilde{\rho}_a+\tilde{w}_a=&\frac{1}{r}\int_{|\lambda_at-r|}^{\lambda_at+r}sg_a(s)ds+\partial_r\int_{-1}^1sm_a(s)\big|_{s=\lambda_at+\theta r}d\theta\\
%  =&\int_{-1}^1sg_a(s)\big|_{s=\lambda_at+\theta r}d\theta+\int_{-1}^1(sm_a)'\big|_{s=\lambda_at+\theta r}\theta d\theta,
%\end{align*}
%then
\begin{align*}
  &\int_0^t\int_{r\leq\varepsilon_0\tau}\Big|(\tilde{\rho}_a+\tilde{w}_a)\cdot r\partial_r(\tilde{\rho}_b+\tilde{w}_b)\Big|drd\tau\\
  \leq&C \int_0^t\int_{\mathbb{R}_+}\Big|\int_{-1}^1sg_a(s)\big|_{s=\lambda_a\tau+\theta r}d\theta+\int_{-1}^1(sm_a)'\big|_{s=\lambda_a\tau+\theta r}\theta d\theta\cdot r(\tilde{\rho}'_{b0}+\tilde{w}'_{b0})\Big|drd\tau\\
\leq &C[\|\xi_\theta g_a\|_{L^1(\mathbb{R}_+)}+\|(\xi_\theta m_a)'\|_{L^1(\mathbb{R}_+)}]\|r\tilde{\rho}'_{b0}+r\tilde{w}'_{b0}\|_{L^1(\mathbb{R}_+)},
\end{align*}
\begin{align*}
  &\int_0^t\int_{r\leq\varepsilon_0\tau}\Big|(\tilde{\rho}_a+\tilde{w}_a)\cdot(\tilde{\rho}_b+\tilde{w}_b)\Big|drd\tau\\
  \leq&C \int_0^t\int_{\mathbb{R}_+}\Big|\int_{-1}^1sg_a(s)\big|_{s=\lambda_a\tau+\theta r}d\theta+\int_{-1}^1(sm_a)'\big|_{s=\lambda_a\tau+\theta r}\theta d\theta\cdot (\tilde{\rho}_{b0}+\tilde{w}_{b0})\Big|drd\tau\\
\leq &C[\|\xi_\theta g_a\|_{L^1(\mathbb{R}_+)}+\|(\xi_\theta m_a)'\|_{L^1(\mathbb{R}_+)}]\|\tilde{\rho}_{b0}+\tilde{w}_{b0}\|_{L^1(\mathbb{R}_+)}.
\end{align*}
And other terms are easier to perform. We have completed the proof of Theorem \ref{thm3}.
\hfill$\Box$
\bigskip

{\bf Acknowledgement} The authors would like to thank Prof. Gui-Qiang Chen for the enlightening discussions and encouragements.

\end{document}